%
%

\input amstex.tex
\input amsppt.sty

\font\footfont=cmr10 at 8  pt 
    
 \font\titfont=cmr10 at 16pt 
 \font\aufont= cmr10 at 13pt


 \magnification=1100
 \NoBlackBoxes
 \nologo





 \define\zz{ {{\bold{Z}_2}}}


 \define\calc{\Cal C}

 \define\calq{\Cal Q}
 
 \define\calz{\Cal Z}



 \define\cycd#1#2{{\calz}_{#1}(#2)}
 \define\cyc#1#2{{\calz}^{#1}(#2)}
 \define\cych#1{{{\calz}^{#1}}}
 \define\cycp#1#2{{\calz}^{#1}(\bbp(#2))}
 \define\cyf#1#2{{\cyc{#1}{#2}}^{fix}}
 \define\cyfd#1#2{{\cycd{#1}{#2}}^{fix}}

 \define\crl#1#2{{\calz}_{\bbr}^{#1}{(#2)}}

 \define\crd#1#2{\widetilde{\calz}_{\bbr}^{#1}{(#2)}}

 \define\crld#1#2{{\calz}_{#1,\bbr}{(#2)}}
 \define\crdd#1#2{\widetilde{\calz}_{#1,{\bbr}}{(#2)}}

 \define\crlh#1{{\calz}_{\bbr}^{#1}}
 \define\crdh#1{{\widetilde{\calz}_{\bbr}^{#1}}}
 \define\cyav#1#2{{{\cyc{#1}{#2}}^{av}}}
 \define\cyavd#1#2{{\cycd{#1}{#2}}^{av}}

 \define\cyaa#1#2{{\cyc{#1}{#2}}^{-}}
 \define\cyaad#1#2{{\cycd{#1}{#2}}^{-}}

 \define\cyq#1#2{{\calq}^{#1}(#2)}
 \define\cyqd#1#2{{\calq}_{#1}(#2)}

 \define\cqt#1#2{{\calz}_{\bbh}^{#1}{(#2)}}
 \define\cqtav#1#2{{\calz}^{#1}{(#2)}^{av}}
 \define\cqtrd#1#2{\widetilde{\calz}_{\bbh}^{#1}{(#2)}}

 \define\cyct#1#2{{\calz}^{#1}(#2)_\zz}
 \define\cyft#1#2{{\cyc{#1}{#2}}^{fix}_\zz}
\define\cxg#1#2{G^{#1}_{\bbc}(#2)}
 \define\reg#1#2{G^{#1}_{\bbr}(#2)}

 \define\cyaat#1#2{{\cyc{#1}{#2}}^{-}_\zz}



 \define\fflag#1#2{{#1}={#1}_{#2} \supset {#1}_{{#2}-1} \supset
 \ldots \supset {#1}_{0} }
 
 \define\vect#1{ {\Cal{V}ect}_{#1}}


 \define\chv#1#2#3{{\calc}^{#1}_{#2}(#3)}
 \define\chvd#1#2#3{{\calc}_{#1,#2}(#3)}
 \define\chm#1#2{{\calc}_{#1}(#2)}



 \define\Claim#1{\subheading{Claim #1}}

 \define\xrightarrow#1{\overset{#1}\to{\rightarrow}}

 \define\cx{\Cal X}


\hfuzz1pc 


\define\bbz{\bold Z}

\define\bbr{\bold R}
\define\bbc{\bold C}
\define\bbh{\bold H}
\define\bbp{\bold P}

\define\ca{\Cal A}
\define\cc{\Cal C}
\define\cd{\Cal D}
\define\ce{\Cal E}
\define\ch{\Cal H}

\define\cl{\Cal L}

\define\cp{\Cal P}
\define\cf{\Cal F}
\define\ck{\Cal K}
\define\cm{\Cal M}

\define\cs{\Cal S}
\define\cz{\Cal Z}
\define\co#1{\Cal O_{#1}}
\define\ct{\Cal T}


\define\a{\alpha}
\redefine\b{\beta}
\define\g{\gamma}
\redefine\d{\delta}

\define\s{\sigma}
\define\z{\zeta}

\define\la{\lambda}
\define\e{\epsilon}

\define\G{\Gamma}

\define\p#1{{\bbp}^{#1}}


\define\blbx{\hfill {\vrule width5pt height5pt depth0pt}}
\redefine\qed{\blbx}

\define\pf{\subheading{Proof}}
\define\Lemma#1{\subheading{Lemma #1}}
\define\Theorem#1{\subheading{Theorem #1}}
\define\Prop#1{\subheading{Proposition #1}}
\define\Cor#1{\subheading{Corollary #1}}
\define\Note#1{\subheading{Note #1}}
\define\Def#1{\subheading{Definition #1}}
\define\Remark#1{\subheading{Remark #1}}
\define\Ex#1{\subheading{Example #1}}
\define\arr{\longrightarrow}

\define\Hom{\operatorname{Hom}}

\redefine\Xi{X_{\infty}}

\define\jac#1#2{\left(\!\!\!\left(
\frac{\partial #1}{\partial #2}
\right)\!\!\!\right)}
\define\restrict#1{\left. #1 \right|_{t_{p+1} = \dots = t_n = 0}}

\define\SP#1#2{{\roman SP}^{#1}(#2)}

\define\coc#1#2#3{\cc^{#1}(#2;\, #3)}
\define\zoc#1#2#3{\cz^{#1}(#2;\, #3)}
\define\zyc#1#2#3{\cz^{#1}(#2 \times #3)}

\define\Div{\roman{ Div}}
\define\ar#1{\overset{#1}\to{\longrightarrow}}

\define\th#1#2{{\Bbb H}^{#1}(#2)}
\define\hth#1#2{\widehat{\Bbb H}^{#1}(#2)}


\define\bad#1#2{\cf_{#2}(#1)}

\define\pch#1{\bbp_{\bbc}(\bbh^{#1})}

\def\z2t{\text{$\bbz_2\ct$}}

\def\hH{\widehat{I\!\!H}}

\def\<{\left<}
\def\>{\right>}
\def\[{\left[}
\def\]{\right]}

\def\wt{\widetilde}
\def\vf{\varphi}

\def\HH{9}

\def\supp{\operatorname{supp \ }}

\redefine\and{\qquad\text{and}\qquad}

\def\cn{\Cal N}

\def\th{^{\text{th}}}

\define\pr{\operatorname{pr}}

\def\wh{\widehat}\

\define\Dp#1{{\cd'}^{#1}(X,p)}
\define\Ep#1{{\ce}^{#1}(X,p)}
\define\Zp#1{{\cz}^{#1}_{\bbz}(X,p)}
\define\Hop#1#2{{H}^{#1}_{#2}(X,p)}
\define\Hp#1{{\hH}^{#1}(X,p)}

\define\Dcx#1{\underline{\bbz}_{\cd}(#1)}
\define\HD#1#2{H^{#1}_{\cd}(X,\bbz(#2))}
       
\redefine\hH{\widehat{\bold H}}

\def\Proj{\operatorname{Proj}}
\def\PH#1{\widehat{#1}}
\def\HH#1#2{\ch_{#1}(#2)}
\def\SH#1{S(#1)}
\def\PHull#1{{\widehat{#1}}_{\operatorname{poly}}}
\def\HPHull#1{{\widehat{#1}}_{\operatorname{hom-poly}}}
\def\tbar{|||}

\def\X{K}
\def\V{\Lambda}
\def\Vaff#1{\Lambda^{0}_{#1}} 
\def\sh{\varphi}
\def\shaff{\psi}  
\def\PSH{\cp\cs\ch} 
\def\SHH{\cs}  
\def\SHHaffine{\cs^0}  


\centerline{\titfont Projective Hulls}\vskip .1in 

\centerline{\titfont and the}\vskip .1in

\centerline{\titfont Projective Gelfand Transform}

\vskip .2in
\centerline{\aufont by}

\vskip .1in
\centerline{\aufont F. Reese Harvey and H. Blaine Lawson, Jr.
\footnote{\footfont Research
 partially supported by the NSF}}
\vskip .3in
\centerline{\sl Dedicated, with affection and deep esteem,}
\smallskip
\centerline{\sl to the memory of  S.-S. Chern.}

\vskip .5in
 
\centerline{\bf Abstract} \medskip
  \font\abstractfont=cmr10 at 10 pt

{{\parindent=.9in\narrower\abstractfont \noindent 
 We introduce the notion of a  
{\bf projective hull}  for subsets of complex projective varieties 
parallel to the idea of a polynomial hull in affine
varieties.  
With this concept, a generalization of J. Wermer's 
classical theorem on the hull of a curve in $\bbc^n$
is established in the projective setting.
The projective hull is shown to have interesting properties and is
related to various extremal functions and capacities in
pluripotential theory. A main analytic result asserts that
for any point $x$ in the projective hull $\PH K$ of a compact set
$K\subset \bbp^n$ there exists a positive current $T$ of bidimension
(1,1) with support in ${\PH K}^-$ and a probability measure $\mu$ on $K$
with $d d^cT=\mu-\delta_x$. This result generalizes to any K\"ahler
manifold and has strong consequences for the structure of $\PH K$.

\ 

\noindent
We also introduce the notion of a {\bf projective spectrum } 
for Banach graded algebras parallel to the Gelfand spectrum of a
Banach algebra. This projective spectrum has universal properties
exactly like those in the Gelfand case. Moreover, the projective hull is
shown to play a role (for graded algebras) completely analogous to that
played by the polynomial hull in the study of finitely generated Banach
algebras.

\ 

\noindent
This paper gives foundations for generalizing many of the results on
boundaries of varieties in $\bbc^n$ to general algebraic manifolds.

}}

\vfill\eject         
               
\centerline{Table of Contents}
\medskip   
       
\hskip .9 in  1.        Introduction
 
\hskip .9 in  2.     The Projective Hull  of a subset of $\bbp^n$     
 
\hskip .9 in  3.      Elementary Properties

\hskip .9 in  4.      The Best Constant  Function,
                 Quasi-plurisubharmonicity  

\hskip .9in    \ \ \ \  and Pluripolarity

\hskip .9 in  5.      The  Picture in Homogeneous Coordinates 
 
\hskip .9 in  6.     The Affine Picture  
 
\hskip .9 in  7.       Theorems of Sadullaev

\hskip .9 in  8.      The Theorem of Fabre

\hskip .9 in  9.    Examples  

\hskip .9 in  \!\! 10.       Compactness and Stability

\hskip .9 in \!\!  11.    Projective Hulls and Jensen Measures (The Main
Theorem)

\hskip .9 in  \!\! 12.     Structure Theorems
  
\hskip .9 in  \!\! 13.  The Projective Spectrum

\hskip .9 in  \!\! 14.  The Projective Gelfand Transform

\hskip .9 in  \!\! 15.  Relation to the Projective Hull
 
\hskip .9 in  \!\! 16.  Finitely generated algebras
 
\hskip .9 in  \!\! 17.  Projective Hulls on Algebraic Manifolds
 
\hskip .9 in  \!\! 18.  Results for General K\"ahler Manifolds
 
\hskip .9 in  \!\! Appendix A.  Norms on $A_*(\bbp^n)$

\vfill\eject

\centerline{\bf 1. Introduction}\medskip

A beautiful classical theorem of John Wermer [W$_1$] states that 
the polynomial hull $\PHull {\g}$ of a compact real analytic curve
$\gamma\subset\bbc^n$, has the property that
 $\PHull {\g} -\g$ is a 1-dimensional
complex analytic subvariety of $\bbc^n-\g$.  (Recall that the polynomial
hull of   $K\subset\subset \bbc^n$ is the set of points $x\in\bbc^n$ such
that $|p(x)|\leq \sup_K|p|$ for all polynomials $p$.)

This paper was largely motivated by
the question: \medskip
\centerline{\sl Does there exist an analogous result for curves in complex
projective space $\bbp^n$?}\medskip

To this end we introduce the notion of the {\sl projective hull}
of a compact set $K\subset\bbp^n$.  It is defined to be the set $\PH K$
of  points $x\in\bbp^n$ for which there exists a constant
$C=C_x$ such that  $$
\|\cp(x)\|\ \leq\ C^d \sup_K\|\cp\|
\tag1.1$$
for all holomorphic sections $\cp$ of $\co{\bbp^n}(d)$  and all $d>0$.
Strong motivation for this definition comes from the fact (Prop. 2.3) 
that  if $\g$ is the boundary of a one-dimensional complex analytic
subvariety $V\subset \bbp^n$, then $V\subseteq \PH {\g}$. Furthermore, for
large classes of non-trivial examples   it is shown in \S 9 that $V =\PH
{\g}$.

The projective hull strictly generalizes the concept of the polynomial
hull in the following sense.  Suppose $K\subset\subset\Omega = $ an
affine open subset of $\bbp^n$.  Then 
$$
{\PH K}_{\operatorname{poly}, \Omega}\ \subseteq \PH K,\ \ \ \ \
\qquad\ \ \  \text{and} \qquad \PH K\subset\subset\Omega
\ \ \Rightarrow\ \ {\PH K}_{\operatorname{poly}, \Omega} = \PH K
$$
\noindent
where ${\PH K}_{\operatorname{poly}, \Omega}$ is defined as above using
the regular functions (polynomials) on $\Omega$. The second statement,
which is non-trivial, is proved in \S 12.
The projective hull also satisfies a Local Maximum Modulus
Principle which states   that for any $K\subset\bbp^n$ and any bounded
domain $U$ in some affine open subset $\Omega$, one has that 
${\PH K} \cap U$ is contained in the $\Omega$-polynomial hull of its
boundary. (See Theorem 12.8 or Theorem 4 below).

The projective hull is always subordinate to the Zariski hull --- if
$K\subset Z\subset \bbp^n$ where $Z$ is an algebraic subvariety, then
$\PH K\subset Z$. Moreover, if  a real curve $\g\subset\bbp^n$ 
is contained in an irreducible algebraic curve $Z$, then $\PH {\g}=Z$.

Note that for  $x\in \PH K$ the infimum of the set of constants $C$ for
which (1.1) holds is again such a constant.  This {\sl best
constant function} $C_K:\PH K\arr\bbr^+$ plays a basic role in the
study of projective hulls. It is bounded iff $\PH K$ is compact, and it
appears repeatedly in many contexts. It is sometimes convenient to extend 
$C_K$ to all of $\bbp^n$  by setting
$C_K(x)=\infty$ for points $x\notin \PH K$.

It is natural to ask for an interpretation of the projective hull
in homogeneous coordinates. Let $\pi:\bbc^{n+1}-\{0\}\arr\bbp^n$ be the
standard projection and for $K\subset \bbp^n$ set $S(K)=\pi^{-1}(K)\cap
S^{2n+1}$.  The polynomial hull of $S(K)$ in $\bbc^{n+1}$ is a compact
subset which is a union of disks centered at the origin. In \S 5 we prove
that
$$
\PH K \ =\ \pi\left\{ \PHull{S(X)} -\{0\}    \right\}
$$
The best constant function $C_K=1/\rho_K$ where $\rho_K(x)$ is the 
radius of the disk in $\PHull{S(X)}$ above $x$.

Interestingly, projective hulls have  already  appeared in a somewhat
hidden way in pluripotential theory. The closest connection is in the work
of  Guedj and Zeriahi [GZ] who (following Demailly)  considered on a
general    K\"ahler manifold $(X,\omega)$ the notion of an
{\sl quasi-plurisubharmonic function}. This is a real-valued function
$v$ on $X$ which satisfies $
dd^c \,v + \omega\ \geq\ 0.
$
The set of these functions is denoted $\PSH_{\omega}(X)$ and for
each compact subset $K\subset X$ there is an associated {\sl extremal
function} 
$$
\Lambda_K(x) \ \equiv \ \sup\left\{v(x): v\in \PSH_{\omega}(X)\ \text{
and } v\bigl|_K\leq 0\right\}.
\tag1.2 $$
Arguments in [GZ] show that for $X-\bbp^n$ the best constant function,
extended to be  $\equiv \infty$ on $\bbp^n-\PH K$, satisfies
$$
\Lambda_K\ =\ \log C_K.
$$

For compact sets $K$ contained in a standard affine coordinate chart
$\bbc^n\subset\bbp^n$ condition (1.1) for $z\in \bbc^n$  is equivalent to
the condition that there exists $C>0$  with
$$
|p(z)|\ \leq\ C^d \sup_K|p|
$$
for all polynomials $p$ of degree $\leq d$ and all $d$.
In this setting the best constant function is related to 
the Siciak extremal function defined in terms of the Lelong class
of subharmonic functions with logarithmic growth [Si]. In particular the
best constant function is finite at exactly those points where
the Siciak function is finite. This is discussed  in \S 6. 

In pluripotential theory it is often customary to regularize extremal
functions to be upper semicontinuous. In the cases of interest here such
regularization gives $\Lambda_K^*\equiv\infty$.  One can
think of our results as showing  that in this situation, the set  where
$\Lambda_K<\infty$ (namely $\PH K$)  has interesting structure, and so
also does $\Lambda_K\bigl|_{\PH K}$.

The condition $\Lambda_K^*\equiv\infty$ is equivalent to $K$
having Bedford-Taylor capacity zero [BT].  It is also equivalent to
$K$ being {\sl pluripolar }, i.e., locally contained in the $-\infty$-set
of a  non-constant plurisubharmonic function (See \S 4). 
This points out the relative subtlety of the projective hull, since there
exist smooth curves in $\bbp^2$ which are not pluripolar [DF].

Another close tie between polynomial and projective hulls comes from the
theory of commutative Banach algebras. In 1941  Gelfand showed that to
every Banach algebra $A$ there is a canonically associated compact
Hausdorff space $X_A$ and a continuous embedding of $A$  into the algebra
$C(X_A)$ of continuous complex-valued functions  on $X_A$. (See [G],    
[Ho] or [AW$_1$].)  The space $X_A$ is universal for  representations of
$A$ in the continuous functions on   compact Hausdorff spaces. The points
of $X_A$ are exactly the representations onto $C(\text{pt})\cong\bbc$,
i.e., the multiplicative linear functionals.

Suppose now that $K$ is a compact subset of $\bbc^n$ and let 
$A(K)$ denote the uniform closure of the polynomials in $C(K)$.
Then there is a canonical homeomorphism
$$
X_{A(K)}\ \cong\ \PHull K
$$
of the Gelfand spectrum with the polynomial hull of $K$.
This engenders a natural correspondence between finitely generated
Banach algebras and polynomially convex subsets of $\bbc^n$, and
enables one to employ the theory of several complex variables in the study
of such algebras.

Now there is a completely parallel story relating projective hulls 
to Banach graded algebras.  This  parallel mimics the relationship
between the Spectrum of a ring and Proj of a graded ring in modern
algebraic geometry. A {\sl Banach graded algebra} is a commutative
$\bbz^+$-graded normed algebra $A_*=\bigoplus_{k\geq 0} A_k$ where each
$A_k$ is a  Banach space. Typical examples are given by:
$A_k=\Gamma_{\text{hol}}(X, \co{}(\la^k))$ with the sup-norm, where $\la $
is a holomorphic hermitian line bundle on a complex manifold $X$, or
$A_k=\Gamma_{\text{cont}}(K, \la^k)$ with the sup-norm,  for a  hermitian
line bundle $\la $  on a compact Hausdorff space $K$. In either case, when
$X=K=$ pt, we have $A_*\cong \bbc[t]$.

For any Banach graded algebra $A_*$ we construct a topological space
$\cx_{A_*}$, called the {\sl projective spectrum} of $A_*$ as the
space of continuous  homomorphisms   $A_*\to \bbc[t]$ divided by the
$\bbc^{\times}$-action corresponding to Aut$(\bbc[t])$.  The space
$\cx_{A_*}$ carries a hermitian line bundle $\la$, and there is a natural
embedding 
$$
A_*\hookrightarrow \bigoplus_{k\geq0}\Gamma(\cx_{A_*}, \la^k)
$$
called the {\sl projective Gelfand transformation}.

Suppose now that $K$ is a compact subset of $\bbp^n$, and let 
$A_*(K)$ denote the restriction to $K$ of the homogeneous coordinate ring
of $\bbp^n$ (represented by homogeneous polynomials in homogeneous
coordinates). Then there is a canonical homeomorphism 
$$
\cx_{A_*(K)}\ \cong\ \PH K
$$
of the projective spectrum with the projective hull of $K$.
This engenders a natural correspondence between finitely generated
Banach graded algebras and projectively convex subsets of $\bbp^n$.

The principal analytic result in this paper is the establishement 
of Jensen measures for points in the projective hull.  The theorem
has a number of interesting corollaries. In particular, with a mild
hypothesis it yields the projective analogue of Wermer's theorem.

Fix $\Lambda>0$ and denote by $\cp_{1,1}(\Lambda)$ the set of positive
currents of bidimension (1,1) with mass $\leq \Lambda$.  For a compact set
$K\subset \bbp^n$ let $\cm_K$ denote the probablitiy measures on $K$ and
let $\PH K(\Lambda)$ denote the set of points in  $x\in \PH K$ with
$\Lambda_K(x)\leq \Lambda$.  \medskip
\noindent
{\bf Main Analytic Theorem.}  {\sl For a compact subset  ${\X}\subset
\bbp^n$ the following are equivalent:

\medskip

\qquad (A)\ \ \ \  $x\in \PH {\X}(\Lambda)$
\medskip

\qquad (B) \ \ There exist  $T\in \cp_{1,1}(\Lambda)$  and
$\mu\in\cm_{\X}$ such that:\medskip 

 \hskip 1.6in (i)    \ \ \ $dd^c T\ =\ \mu - \d_x$
\smallskip

 \hskip 1.6in (ii)    \ \ $\supp(T)\ \subset {\PH {\X}}^{-} \equiv
\text{the closure of } \PH \X$  \smallskip }

\medskip

\Note{}  This theorem was inspired by a result of
Duval-Sibony [DS, Thm. 4.2] in the affine case, and our proof
incorporates their Hahn-Banach technniques. However, much more is
required.  Our projective result is (necessarily) quantitative.
Furthermore, one must work in this case to find a current $T$ with
support in the closure of the projective hull.

It is tempting to apply   [DS] directly by considering the set
$S(\X)\subset S^{2n+1}\subset\bbc^{n+1}$ in homogeneous coordinates
discussed above, and then push their positive (1,1)-current $T$ forward to
$\bbp^n$ by the projection.  However, there is nothing in [DS] that
indicates how to construct $T$ so that $0\notin \supp(T)$. Indeed in
homogeneous coordinates, much of the subtlety of this subject takes place
near the origin.

\medskip
As a corollary of the Main Analytic Result one can show that  for any 
$x\in {\PH {\X}}^-$  there  are probability measures $\nu\in \cm_{{\PH
{\X}}^-}$, $\mu\in \cm_{\X}$ and a current  $T\in\cp_{1,1}$  
 with
 $dd^c T\ =\ \mu - \nu$ and 
  $x\in \supp(T)\ \subset {\PH {\X}}^{-}$.

Note that a current $T_x\in \cp_{1,1}$ of least mass with support in  
${\PH {\X}}^{-}$ and $dd^c T_x\ =\ \mu - \d_x$ satisfies
$$
M(T_x)\ =\ \Lambda_K(x)
$$

\medskip
Fix a compact subset $K\subset \bbp^n$ . One of the important consequences
of the main theorem is the following.

\Theorem{1}{\sl The set 
${\PH K}^- -K$ is 1-concave in $\bbp^n-K$.}
\medskip

One-concavity is a strong local condition which means essentially no
local peak points under holomorphic maps to $\bbc$.  (The definition 
is given in \S 12.) It has the following immediate consequence.

\Cor{1}{\sl 
${\PH K}^- -K$ has locally positive Hausdorff 2-measure.}
\medskip
 
Using Theorem 1 combined  with work of Dinh and Lawrence [DL] or Sibony
[Sib, Thm. 17]   we conclude the following.

\Theorem{2}{\sl If the  Hausdorff 2-measure of ${\PH K}^- -K$ is locally
finite on some open subset $U\subset \bbp^n-K$, then $({\PH K}^- -K)\cap U$
is a 1-dimensional complex analytic subvariety of $U$.}
\medskip

\medskip

\Theorem {3}{\sl Suppose that 
${\X}\subset \subset\bbc^n$ and let $\PH {\X}_0$
be a connected component of ${\PH {\X}}^- -{\X}$  which is bounded in
$\bbc^n$.  Then $\PH {\X}_0$ is contained in the polynomial hull of ${\X}$.}

\Cor{3} {\sl
If $\PH {\X}\subset \subset \Omega=\bbp^n-D$ for some algebraic
hypersurface $D$, then $\PH {\X} ={\PH {\X}}_{\Omega}$.  In particular if
$\g\subset\Omega$ is a C$^1$-curve and if $\PH \g\subset\subset \Omega$,
then $\PH \g -\g$  is a 1-dimensional analytic subvariety of $\Omega-\g$.
}

\medskip

One might conjecture that if $\PH \g$ is not an algebraic subvariety, then
it is contained in the complement of some divisor, and Corollary 2 would
give a projective version of Wermer's Theorem.  However, recent beautiful
work of Bruno Fabre [Fa$_1$] shows that this is far from true. His results
are discussed in \S 8.

Using our Main Theorem we establish the following  local structure theorem
which yields, in particular, the Local Maximum Modulus Property for
projective hulls.

\Theorem{4}{\sl For any
bounded domain $U\subset\subset\Omega^{\text{affine open}}
\subset \bbp^n$,} 
$$
{\PH K}^- \cap U\ \subseteq\ 
{\biggl\{(\PH K \cap \partial U) \  \cup\ (K\cap
U)\biggl\}}_{\operatorname{poly}, \Omega}^{\widehat { \ \ \ \  }} $$

\medskip

We also obtain the following generalization of Wermer's Theorem.

\Theorem{5} {\sl Let $\gamma\subset\bbp^n$ be a finite union of real
analytic curves. Then $\PH{\gamma}$ has Hausdorff dimension 2.
Furthermore, if the Hausdorff 2-measure of $\PH{\gamma}^-$ is finite in a
neighborhood of some algebraic hypersurface, then $\PH {\g}-\g$ is a
1-dimensional complex analytic subvariety of $\bbp^n-\g$.

The same conlcusion holds for any smooth pluripolar curve $\g$  in
$\bbp^2$.}
\medskip

Theorems 1--5 are proved in \S 12.

Given a complex manifold $X$  and  a hermitian holomorphic line bundle
$\la\to X$ there is an analogue  of the projective hull of $K\subset\subset
X$ defined to be the set ${\PH K}_X$ of points $x\in X$ satisfying (1.1)
for all $\cp\in  H^0(X, \co{}(\la^d))$ and all $d>0$. There is also an
analogue $\Lambda_{K,X}$ of the extremal function (1.2). This is
discussed in  \S 17 where we prove that {\bf the projective hull is 
 intrinsic}, namely:

\Theorem{6} {\sl Suppose that  $X\subset \bbp^N$ is an algebraic manifold
and  $\la=\co{X}(1)$.  Then for any compact subset $K\subset X$ we have
$$
\Lambda_{K,X} \ =\ \Lambda_K\bigl|_X
\and
{\PH K}_X\ =\ \PH K.
$$
}

The reader may recall that   extremal functions can be defined on any
K\"ahler manifold $(X, \omega)$ using the quasi-plurisubharmonic
functions by (1.2). One can therefore define the {\bf  $\omega$-hull}
$\PH K$ of a compact subset $K\subset X$ to be the set of points $x\in X$
where  $\Lambda_K(x)<\infty$.

In \S 18 we establish   the following.

\medskip
\noindent
{\bf  Main Analytic Theorem for Arbitrary
K\"ahler Manifolds  } {\sl   Let $X$ be any K\"ahler manifold.      
Then for any compact subset $K\subset X$ with  $\PH K\subset\subset X$ the
following are equivalent.  \medskip 

\qquad (A)\ \ \ $x \ \in\ \PH K(\Lambda)$\medskip

\qquad (B) \ \ There exist $T\in \cp_{1,1}(X)$ with $M(T)\leq \Lambda$ and
$\supp(T)\subset {\PH K}^-$ such that
$$
dd^c T \ =\ \mu - \delta_x
$$
\qquad\qquad \ \ \ \ \ where $\mu \in \cm_K$.
}
 
\medskip

The proof of this result is more
rounded and conceptual than the one given for the special case in \S 11. 
Most of the consequences of the special case cited above  carry over to
the general setting.

We point out that this paper lays the foundation for a new
characterization of boundaries of subvarieties in a compact Kahler
manifold $X$. For $X=\bbp^n$ this problem has been studied in [Do],
[DH$_{1,2}$], and [HL$_3$] where such boundaries were characterized in
terms of  analytic transforms and non-linear moment conditions.  However,
using   results  in this paper, the authors have formulated a quite
different characterization of the boundaries of {\sl positive} holomorphic
chains in terms of projective linking numbers [HL$_4$].
This generalizes the work of Alexander and Wermer [AW$_2$], [W$_2$] in
$\bbc^n$. The results in [HL$_4$] also cover the case of a general
K\"ahler manifold $X$.

The authors would like to particularly thank Eric Bedford 
for several very useful conversations relating to this article.  We also
thank Tien-Chong Dinh and Vincent Guedj for explaining their results which
have played an important role here. We are indebted to Nessim Sibony, 
Vincent Guedj, Ahmed Zeriahi, and Tien-Cuong Dinh for numerous useful
comments.  The second author would like to thank the Institut Henri
Poincar\'e and in particular Gennadi Henkin and Nessim Sibony for their
hospitality during the development of this work.

\vskip .3in


\centerline{\bf 2. The Projective Hull   of a subset of $\bbp^n$}\medskip

Let $\co {}(d) \arr \bbp^n$ denote the holomorphic line bundle of Chern
class $d$ over complex projective $n$-space. Note that any hermitian
metric on $\co{}(-1)$ naturally  induces a hermitian metric on $\co{}(d)$
for each $d\in\bbz$.  This family of metrics has the property that 
$$
 \qquad\qquad\qquad\qquad\qquad\qquad
\|v^{\otimes d}\|\ =\ \|v\|^d\qquad\qquad \text{ for any } v\in \co{}(d_0),
\text{ any } d_0
\tag2.1$$
and $|(v,w)|=\|v\|\cdot\|w\|$ for $v\in \co{}(d), w\in \co{}(-d)$
for any $d$.  Fix any such family of metrics and consider a compact subset
${\X}\subset \bbp^n$.

\Def{2.1}  The {\bf projective hull of } ${\X}$ is the subset $\PH {\X}$ of
points   $x\in \bbp^n$ with the following property:
There exists a constant $C$  (depending on $x$)  such that
$$
\|\sigma(x)\|\ \leq \ C^d\sup_{\X}\|\sigma\|
\tag2.2$$
for all $\s \in H^0(\bbp^n;\co{}(d))$ and all $d\geq 0$.

The infimum of all $C$ for which (2.2) holds will be called the {\sl best
constant function} and denoted by $C_{\X}(x)$.

\medskip

\noindent
{\bf Exercise 2.2} The   hull $\PH {\X}$ is independent of the
choice of hermitian metric on $\co{}(-1)$. \medskip

The following fact was a primary motivation for considering this
concept.

\Prop{2.3} {\sl Let $V$ be a compact connected Riemann surface with
boundary $dV\neq \emptyset$. Suppose $f:V\to \bbp^n$ is a holomorphic map
which extends holomorphically across the boundary.
Then 
$$ f(V)\
\subseteq \PH {f(dV)} 
\tag2.3$$
}

\pf By assumption  $V\subset\subset \wt V$  for
some connected non-compact Riemann surface $\wt V$ and $f$ extends
holomorphically to $\wt V$. 
Since $\wt V$ is Stein, there is a holomorphic trivialization
of the pull-back
$
f^*\co{}(1)\  @>{\cong}>>\ \wt V\times \bbc
$
which yields trivializations
$$
f^*\co{}(d)\  @>{\cong}>>\ \wt V\times \bbc.
\tag2.4
$$
for all $d\geq 1$.
With respect to (2.4) the pull-back metric is of the form
$$
\|v\| = \lambda^d |v| 
$$
for some smooth function  $\lambda:\wt V\to \bbr^+$.

Fix   $p\in V-dV$ and $\sigma \in H^0(\bbp^n,\co{}(d))$.  Let $\wt
\s:\wt V\to \bbc$ be the holomorphic function corresponding to 
$f^*\s$ under (2.4). Then applying the maximum principle to 
the compact subdomain $V\subset \wt V$ gives
$$
\|\wt \s (p)\|=\lambda(p)^d|\wt \s(p)|\leq
\lambda(p)^d\sup_{dV}|\wt\s| \leq 
\left(\frac{\lambda(p)}{\mu}\right)^d\sup_{dV}\|\wt\s\|,
\qquad
\text{ where }
\mu \equiv \inf_{dV}\lambda
$$
as desired. \qed

This proof actually establishes the following.

\Prop{2.4} {\sl The conclusion (2.3) holds for any  map
$f:V\to \bbp^n$,  holomorphic on $V-dV$ and continuous on $V$, such that
the pull-back $f^*\co{}(1)$ admits a trivialization which is holomorphic
on $V-dV$ and continuous on $V$.

}

\Remark{2.5} Although the projective hull is independent of the metric
chosen on  $\co{}(-1)$ it is convenient to work with the  standard
metric defined as follows.  Recall that 
$$
\co{}(-1) = \{(\ell, v)\in \bbp^n\times \bbc^{n+1}\ :\ v\in \ell\}
\tag2.5
$$
and projection to the second factor gives a map $\pr_2:\co{}(-1)\to
\bbc^{n+1}$ which is an isomorphism outside the zero-section and collapses
the zero-section to the origin (the blow-up of 0). The {\bf standard
metric on $\co{}(-1)$} is the unique hermitian metric whose unit circle
bundle corresponds to the unit sphere $S^{2n+1}\subset \bbc^{n+1}$ under
the map $\pr_2$.

Recall that any section $\s\in H^0(\bbp^n,\co{}(d))$ gives a function
$\wt \s :\co{}(-1)\to \bbc$ which of degree $d$ on each fibre and descends
under $\pr_2$ to a homogeneous polynmial $\wt\s:\bbc^{n+1}\to \bbc$ of
degree $d$. This gives the identification
$$\qquad 
H^0(\bbp^n,\co{}(d)) \ \cong\ \bbc[Z_0,...,Z_n]_d, \qquad \ \  d\geq 0
\tag2.6
$$
with the space of homogeneous polynomials of degree $d$ in homogeneous
coordinates.

Given a polynomial $P\in \bbc[Z_0,...,Z_n]_d$, its standard norm at 
$x=[Z]\in \bbp^n$, when considered as a section $\cp$  of $\co{}(d)$, is 
$$
\|\cp(x)\| \ =\ \frac{|P(Z)|}{\|Z\|^d}
\tag2.7$$
In particular, given a subset ${\X}\subset \bbp^n$, let $\SH {\X}
=\{z\in S^{2n+1} : \pi(z) \in {\X}\}$ where $\pi:\bbc^{n+1}-\{0\}\to \bbp^n$
is the homogeneous coordinate map.  Then
$$
\sup_{\X} \|\cp\|\ =\ \sup_{\SH {\X}} |P|
\tag2.8$$
for $P$, $\cp$ as above.

\vfill\eject


\centerline{\bf 3. Elementary Properties}\medskip
 The projective hull has several nice features.

\Prop{3.1} {\sl

\qquad(i)\ \  \ If ${\X}\subseteq Y$, then $\PH {\X}\subseteq \PH Y$.

\qquad(ii)\ \  If $Y$ is an algebraic subvariety, then $\PH Y=Y$.

\qquad(iii)\ \   For any ${\X}\subset \bbp^n$, $\PH {\X}$ is contained in
the Zariski hull of ${\X}$.  }

\pf Part (i) is clear and (ii) $\Rightarrow $ (iii). To prove (ii) it 
suffices to show that if $D\subset \bbp^n$ is an algebraic hypersurface
and $Y\subset D$, then $\PH Y\subset D$. Write $D$
as $D=$ Div$(\s)$ for some  $\s\in H^0(\bbp^n,\co{}(d))$.
Then, $Y\subset D \ \ \Leftrightarrow \ \  \s\bigr|_Y =0  \ \ \Rightarrow
\ \  \s\bigr|_{\PH Y} =0  \ \ \Leftrightarrow \ \  \PH Y\subset D$.\qed 

\medskip
The next result says that taking projective hulls commutes with
Veronese re-embeddings.

\Prop{3.2} {\sl Let ${\PH {\X}}_k$ denote the set of points $x\in \bbp^n$
for which there exists $C=C(x)$ such that 
$\|\s(x)\|\leq C^d\sup_{\X}\|\s\|$
for all $\s\in H^0(\bbp^n,\co{}(dk))$ and all $d\geq1$. Then}
$$
{\PH {\X}}_k\ =\ \PH {\X}
$$
\pf 
Suppose $x\in \PH {\X}$ and let $C=C(x)$ be the constant given in
Definition 2.1. Then $\|\s(x)\|\leq (C^k)^d\sup_{\X}\|\s\|$ for all $\s\in
H^0(\bbp^n,\co{}(dk))$ and so $x\in {\PH {\X}}_k$.

On the other hand, suppose $x\in {\PH {\X}}_k$ and let $C=C(x)$ be the
constant in the definition above. Suppose $\s\in H^0(\bbp^n,\co{}(d))$
is given. Then $\s(x)^{\otimes k}\in H^0(\bbp^n,\co{}(dk))$ and so
$$
\|\s(x)\|^{ k} = \|\s(x)^{\otimes k}\| \leq C^d\sup_{\X}\|\s^{\otimes k}\|
=C^d\sup_{\X}\|\s\|^k=C_0^{dk}\left\{\sup_{\X}\|\s\|\right\}^k
$$  
where $C_0\equiv C^{1/k}$. Taking $k$th roots shows that $x\in \PH {\X}$.\qed
\medskip

Recall that for any complex manifold $\Omega $ and any subset
${\X}\subset\Omega$, the {\bf  holomorphic hull of
${\X}$ in $\Omega$ } is the set
$$
\HH{\Omega}{\X}\ \equiv \ \{x\in \Omega\  :\  |f(x)|\leq \sup_{\X}|f|
\ \ \text{ for all }f\in \co {\X}\}
$$

\Prop{3.3}  {\sl Let $\Omega \subset \bbp^n$ be any open subset containing
${\X}$ with the property that Image$\{H^1(\bbp^n,\co{}^{\times})\to
H^1(\Omega,\co{}^{\times})\}$ is finite. Then  }
$$
\HH{\Omega}{\X}\ \subset \PH {\X}.
$$

\Cor{3.4}  
$$
\PH {\X}\ \supseteq\ \bigcup_{\Omega=\bbp-D}\HH{\Omega}{\X}
$$
{\sl where $D$ ranges over all divisors in $\bbp^n$}.
\medskip
\noindent
{\bf Proof that 3.3 $\Rightarrow$ 3.4.}   One has that
 Image$\{H^1(\bbp^n,\co{}^{\times})\to
H^1(\bbp^n-D,\co{}^{\times})\}  \cong \bbz/k$ where $k=$ degree$(D)$.\qed

\medskip
\noindent
{\bf Proof of 3.3.} By assumption there is an integer $k>0$ such that 
$\co{}(k)\bigr|_{\Omega}$ is trivial. The argument given for Proposition
2.3 applies directly to prove that $\HH \Omega {\X} \subseteq
{\PH {\X}}_k\ =\ \PH {\X}$

\Remark{3.5} It should be noted that the containment in Corollary 3.4 is
{\bf not an equality} in general, even if one assumes that ${\X}$ is a
smoothly embedded $S^1$. See \S 8.

\vfill\eject

\centerline{\bf 4. The Best Constant  Function,
                 Quasi-plurisubharmonicity  and Pluripolarity}\medskip

Definition 2.1 leads  naturally  to considering the family
$\SHH_{\X}$ of functions:
$$
\sh = \frac 1d \log \|\cp\|, \qquad  \cp \in H^0(\bbp^n, \co{}(d))\
\text{\ for\ } \ d > 0\ 
\tag4.1$$
with the property that 
$$
\sh\ \leq \ 0 \qquad\text{ on }\ {\X}.
\tag4.2
$$
The  {\bf associated extremal function} 
$$
\V_{\X}(x)\ \equiv\ \sup_{\sh\in \SHH_{\X}}  \sh(x)
$$ 
is the log of the {\bf best constant 
$$
C_{\X}(x)\ =\ \exp(\V_{\X}(x))
\tag4.3$$
satisying (2.2)}. In particular
$$
\PH {\X}\ =\ \left\{x\in\bbp^n: \V_{\X}(x)
<\infty    \right\}
\tag4.4$$

Each function $\sh=\frac 1d\log\|\cp\| \in \SHH_{\X}$  satisfies the
current equation:
$$
d d^c \sh\ =\ \Div(\cp)-\omega \qquad \text{ on }\bbp^n
$$
where the (1,1)-form $\omega$ is the standard K\"ahler form on $\bbp^n$.
These important functions sit in the following, much larger, convex
cone introduced by Demailly.
\Def{4.1}  A function $\sh\in L^1(\bbp^n)$ is called 
{\bf quasi-plurisubharmonic} (or {\bf $\omega$-quasi-plurisubharmonic})
if it is an upper-semicontinuous, $[-\infty, \infty)$-valued function
which satisfies
$$
d d^c \sh+\omega \ \geq \ 0 \qquad \text{ on }\bbp^n.
\tag4.5
$$
The set of these functions will be denoted by ${\PSH}_{{\omega}}$.
\medskip

Note that ${\PSH}_{{\omega}}$ contains $C^\infty$-functions as well as the
highly singular ones in (4.1).

The following useful result can be found in [GZ, Proof of Thm. 4.2].
Analogues for general K\"ahler manifolds follow from work of Demailly
[D$_*$].

\Prop {4.2} {\sl For any compact subset $K\subset \bbp^n$
$$
{\V}_{\X}(x)\ \equiv\ \sup \ \left
\{ \sh(x) : \ \sh\in {\PSH}_{\omega}\text{\  and \ }\sh\bigl|_K\ \leq\
0\right\} $$
}

Recall that a  Borel subset ${\X}\subset \bbp^n$ is called {\bf  globally 
$\omega$-pluripolar} if ${\X}\subseteq\{x\in \bbp^n: \sh(x)=-\infty\}$ for
some  quasi-plurisubharmonic function $\sh \leq 0$ which is not
identically $-\infty$ on $\bbp^n$. The set ${\X}$ is called (locally)
{\bf  pluripolar} if every point $x\in {\X}$ has a connected neighborhood
$\co{}$ such that ${\X}\cap \co{}\subseteq\{x\in \bbp^n: \sh(x)=-\infty\}$
for some plurisubharmonic function $\sh$ on $\co{}$ which is not
identically $-\infty$. 

Guedj and Zeriahi  introduced
and systematically studied  quasi-plurisubharmonic functions
in [GZ]. They also considered a notion of  $\omega$-capacity, due
originally to Dinh-Sibony [DiS], and related it to capacities of
Bedford-Taylor [BT], Alexander [A$_2$],  Sibony-Wong [SW] and others. In
all cases the sets of capacity zero are shown to be the same, and the
following is proved.

\Theorem {4.3. (Guedj-Zeriahi [GZ])} {\sl Let ${\X}\subset \bbp^n$ be  a
compact subset, and denote by ${\V}_{\X}^*$ the upper-semicontinuous
regularization of the function ${\V}_{\X}$.   Then the following are
equivalent:

\qquad \qquad \qquad \qquad (1)\ \ \ $\sup_{\,\bbp^n} {\V}_{\X}\ =\
\infty$. \smallskip

\qquad \qquad \qquad\qquad (2)\ \ \ ${\V}_{\X}^*\ \equiv \-\infty$.
\smallskip

\qquad \qquad\qquad \qquad (3)\ \ \ ${\X}$ has capacity zero.\smallskip

\qquad \qquad\qquad \qquad (4)\ \ \ ${\X}$ is globally
$\omega$-pluripolar.\smallskip

\qquad \qquad\qquad \qquad (5)\ \ \ ${\X}$ is pluripolar.\medskip
\noindent
 }

\Cor{4.4} {\sl For a compact subset ${\X}\subset \bbp^n$ \medskip
$$
{\X} \text{ is pluripolar}\qquad\Leftrightarrow\qquad
\PH {\X} \text{ is pluripolar}\qquad\Leftrightarrow\qquad
\PH {\X}\neq \bbp^n
$$
}
\pf
If ${\X}$ is pluripolar, then by Theorem 4.3(4) ${\X}$ is contained in the
$-\infty$ locus of a negative quasi-plurisubharmonic function $\sh\not
\equiv -\infty$ on $\bbp^n$. Now by (4.3) and Proposition 4.2, $x\in \PH
{\X}$ iff $\sh(x)\leq \sup_{\X} \sh + {\V}_{\X}(x)$ for all 
quasi-plurisubharmonic functions $\sh$ on $\bbp^n$. Hence, $\sh\equiv
-\infty$ on $\PH {\X}$ and so $\PH {\X}$ is $\omega$-pluripolar.  The
converse is obvious, and the first equivalence is established.

Evidently if $\PH {\X}$ is pluripolar, then $\PH {\X} \neq \bbp^n$.
However, if $\PH {\X}\neq \bbp^n$, then by definition ${\V}_{\X} $ is
unbounded on $\bbp^n$ and hence ${\X}$ is pluripolar by Theorem 4.3.
\qed

\Remark{4.5} This corollary highlights the delicate nature of the
projective hull. It is known (cf. [DF],[CLP]) that: 
$$
\text{There exist $C^{\infty}$ closed curves
in $\bbp^2$ which  are not pluripolar.}
$$
The example in [CLP] actually bounds a holomorphic disk in $\bbc^2$.
For this curve, the polynomial hull is the holomorphic disk and
the projective hull is all of $\bbp^2$.

Note however that any real analytic curve is pluripolar, and therefore
its projective hull is a proper subset of $\bbp^n$. A nice
characterization of smooth graphs over the circle which are
pluripolar is given in [CLP].

\Note{4.6} A version of Theorem 4.3 is established in [GZ] with $\bbp^n$
replaced by any compact K\"ahler manifold (cf. \S 18).

\Def{4.7} By the {\bf pluripolar hull} of a Borel set $K\subset \bbp^n$
we mean the set $\PH K_{\text{pp}}$ of points $x\in \bbp^n$ with the
property that $\varphi(x)=-\infty$ for every non-constant  $\varphi \in
\PSH_{\omega}$ with $\varphi\bigl|_K\equiv -\infty$.

\Prop{4.8} {\sl For any compact subset $K\subset\bbp^n$, one has 
$\PH K \subseteq \PH K_{\text{pp}}$.}
\pf Suppose $K\subset \varphi^{-1}(-\infty)$ for a non-constant
$\varphi\in\PSH_{\omega}$. Then for every $c\in\bbr$, $\varphi+c\leq 0$ on
$K$, and so $\varphi(x)+c\leq \Lambda_K(x) < \infty$. Hence, $\varphi(x)
\leq \Lambda_K(x)-c$ for all $c\in\bbr$.\qed


\vfill\eject

\centerline{\bf 5. The  Picture in Homogeneous Coordinates.}\medskip
Consider homogeneous cooordinates $\pi:\bbc^{n+1}-\{0\}\to \bbp^n$
and endow $\bbc^{n+1}$ with the standard hermitian metric.
For any subset ${\X}\subset \bbp^n$ set
$$
\SH {\X}\ \equiv \ \pi^{-1}({\X})\cap S^{2n+1}
$$
where $S^{2n+1}\subset \bbc^{n+1}$ is the unit sphere. 
Note that $\SH {\X}$ is $S^1$-invariant, where $S^1\subset \bbc$ acts by
scalar multiplication. 
In this section we shall characterize the polynomial hull
$\PH {\X}$  in terms of $\SH {\X}$.

Recall   (2.7) that: $\|\cp(x)\|=|P(Z)|/\|Z\|^d$ if $Z\in \bbc^{n+1}-\{0\}$
and $x=\pi Z$, where $\cp$ is the section of $\co{}(d)$ corresponding
to the homogeneous polynomial $P\in \bbc[Z]_d$ of degree $d$.  Fix
$x\in\bbp^n$.  By definition $x\in \PH {\X}$ if and only if
$$
\|\cp(x)\|\ \leq \ C^d\sup_{\X}\|\cp\|
\tag5.1
$$
for all $\cp\in H^0(\bbp^n, \cp{}(d))$ and $d\geq0$.
This condition can be restated in homogeneous coordinates as
$$
|P(Z)|\ \leq\ \sup_{S({\X})}|P| 
\tag5.2
$$
for all $ P\in\bbc[Z]_d, d\geq0$ and for all $Z\in \bbc^{n+1}$ with
$\pi(Z)=x$ and $\|Z\|\leq 1/C$.

To see that (5.1) and (5.2) are equivalent recall that $\sup_{\X}\|\cp\|
=\sup_{S({\X})}|P|$ by (2.8). Substituting into (5.1) yields
$$
\frac{|P(Z)|}{\|Z\|^d}\ \leq \ C^d\sup_{S({\X})}|P|.
\tag5.3
$$
Hence, (5.1) implies (5.2). Now (5.2), with $Z\in \pi^{-1}x$ chosen so
that $\|Z\|=1/C$, is exactly (5.3), which implies (5.1)

\Def{5.1} The {\bf homogeneous polynomial hull} of a subset $K\subset
\bbc^{n+1}$ is the set $\HPHull K$ of points $Z\in \bbc^{n+1}$ with the
property that 
$$
 |P(Z)|\leq \sup_K |P|  
$$
for all homogeneous polynomials $P$.
\medskip

Given ${\X}\subset \bbp^n$,
let $\rho:\PH {\X} \to \bbr^+$ be defined by 
$$
\rho(x)\ \equiv\ \sup\left\{\|Z\| : \pi(Z)=x\ \ \text{ and }\ \ Z\in
\HPHull{S({\X})}\right\}. \tag5.4
$$
that is, the radius of the largest disk about zero in the line
$\pi^{-1}(x)$ which is contained in the homogeneous polynomial hull of
$S({\X})$.

\Prop{5.2} {\sl For ${\X}\subset \bbp^n$
$$
\rho(x)\ =\ \frac{1}{C_{\X}(x)}\qquad\text{for all  }\ x\in \PH {\X}
$$
where $C_{\X}$ is the best constant function (cf. (2.1)).
}
\pf  This is immediate from the equivalence of conditions (5.1) and (5.2)
above.\qed

\Cor{5.3}  {\sl For any subset ${\X}\subset\bbp^n$}
$$
\PH {\X} \ =\ \pi\left\{\HPHull {S({\X})}  - \{0\}  \right\}
$$

The following result is classical (cf. [A$_2$]). We include a proof for
completeness.

\Prop{5.4}  {\sl For any $S^1$-invariant subset $K\subset\bbc^{n+1}$ 
$$
\HPHull K \ =\ \PHull K
$$
where $\PHull K$ denotes the ordinary polynomial hull of $K$.}

\pf
Clearly $\PHull K\subseteq \HPHull K$, and we need only prove 
$\HPHull K\subseteq \PHull K$. For this we use the following Lemma.

\Lemma{5.5} {\sl Let $P$ be a polynomial in $\bbc^{n+1}$ and write
$P=\sum_{m=0}^N P_m$ where $P_m$ is homogeneous of degree $m$. Then for 
any $Z\in \bbc^{n+1}$}
$$
|P_m(Z)|\ \leq\ \sup_{\theta}\left|P(e^{i\theta}Z)\right|
$$
\pf
Note that 
$
P(\lambda Z)= \sum_{m=0}^N\lambda^m P_m(Z),
$
for   $\lambda\in\bbc$, and therefore
$$
P_m(Z)\ =\ 
\frac1{m!}\frac{\partial^m}{\partial \lambda^m} P(\lambda Z)\biggr|_{\la=0}
\ =\ \frac1{2\pi i}\int_{|\la|=1} \frac{P(\la Z)\,d\la}{\la^{m+1}}
$$
  which gives that
$$
|P_m(Z)|\ \leq\ 
\frac1{2\pi }\int_{0}^{2\pi}\, \left|P(e^{i\theta}Z)\right| \, d\theta 
\ \leq\ \sup_{\theta}\left|P(e^{i\theta}Z)\right|.  
$$
\qed

\Cor{5.6}  {\sl Let  $P=\sum_{m=0}^N P_m$ be as in Lemma 5.5, and suppose
$K\subset \bbc^{n+1}$ is $S^1$-invariant.  Then for $0\leq m\leq N$ one has}
$$
\sup_K |P_m|\ \leq\ \sup_K |P|
$$

Suppose now that $Z\in \HPHull K$ and let  $P=\sum_{m=0}^N P_m$ be any
polynomial decomposed as above.  Then 
$$
|P(Z)|\ \leq \sum_{m=0}^N|P_m(Z)|\ \leq \ \sum_{m=0}^N\sup_K |P_m|
\ \leq \ \sum_{m=0}^N\sup_K |P|\ =\ (N+1)\sup_K|P|.
$$
Applying this to $P^q$ gives $$|P(Z)|^q\leq (N+1)q \sup_K|P|^q=
 (N+1)q (\sup_K|P|)^q$$
and therefore
$$
|P(Z)|\ \leq\ (Nq+q)^{\frac1q} \sup_K|P|.
$$
Since $\lim_{q\to\infty}(Nq+q)^{\frac1q}=1$, we have
$|P(Z)|\ \leq\  \sup_K|P|$, and so $Z\in \PHull K$. \qed
\medskip

Combining Corollary 5.3 and Proposition 5.4 gives the following.

\Theorem{5.7} {\sl For any subset ${\X}\subset \bbp^n$ let $\PH {\X}\subset
\bbp^n$ denote its projective hull.  Then}
$$
\PH {\X} \ =\ \pi\left\{\PHull {S({\X})} - \{0\} \right\}
$$

\vskip .3in

\centerline{\bf 6. The Affine Picture.}\medskip
Let ${\X} \subset \bbp^n$ be a compact subset contained in an affine open
chart $\bbc^n=\bbp^n-\bbp^{n-1}$.
\Prop{6.1}{\sl A point $z\in \bbc^n$ lies in $\PH {\X}$ if and only if 
there exists a constant $c>0$ such that
$$
|p(z)|\ \leq\ c^{d} \sup_{\X} |p|
\tag6.1
$$
for all polynomials $p\in \bbc[z_1,...,z_n]$ of degree $\leq d$.
}

\pf Choose homogeneous coordinates $[Z_0:\dots:Z_n]$ for $\bbp^n$
such that $\bbp^{n-1}=\{Z_0=0\}$.  Choosing affine coordinates 
$(z_1,...,z_n)\mapsto[1:z_1:\dots:z_n]$ we identify $H^0(\bbp^n,\co{}(d))$
with the space $\bbc[z_1,...,z_n]_{\leq d}$ of polynomials of degree $\leq
d$ (cf. (2.6)). (This results from the trivialization of $\co{}(d)$ over
$\bbc^n$ via the section $Z_0^d$.)\  In this picture the standard metric
has the form
$$
\|\cp\|_z\ =\ \frac{|p(z)|}{(1+\|z\|^2)^{\frac d2}}.
\tag6.2
$$ 
where $\cp$ denotes the section corresponding to $p$. Consequently a
point $z\in \bbc^n$ lies in $\PH {\X}$ iff there is a constant $C$ such
that  $$
\frac{|p(z)|}{(1+\|z\|^2)^{\frac d2}}\ \leq\ C^d\sup_{\zeta\in {\X}}
\left\{\frac{|p(\zeta)|}{(1+\|\zeta\|^2)^{\frac d2}}\right\}
\tag6.3
$$
from which the result follows directly.\qed
\medskip

In affine coordinates the family of functions $\SHH_{\X}$ defined by 
(4.1) and (4.2) is given by
$$
\sh\ =\ {\tsize\frac 1d}\log |p(z)|-\log\sqrt{1+|z|^2}
\tag6.4$$
where $p\in\bbc[z_1,...,z_n]_{\leq d}$, with the property that
$$
\sh\ \leq \ 0\qquad\text{on}\ \ {\X}\subset \bbc^n.
\tag6.5$$

Proposition 6.1 suggests we consider the family $\SHHaffine_{\X}$ of
plurisubharmonic functions
$$
\shaff(z)\ =\ {\tsize\frac 1d}\log |p(z)|
\tag6.6$$
for $p\in \bbc[z_1,...,z_n]_{\leq d}$, with the property that
$$
\shaff\ \leq \ 0\qquad\text{on}\ \ {\X},
\tag6.7$$
and the associated extremal function
$$
{\Vaff{\X}}(z)\ \equiv\ \sup_{\shaff\in\SHHaffine_{\X}} \shaff(z) 
$$
on $\bbc^n$. Note that for $x\in \bbc^n$, 
$$
\Vaff{\X}(x)\ <\ \infty\qquad\text{ iff }\qquad {\V}_{\X}(x)\ <\ \infty
\tag6.8$$
Therefore, from (4.4) we have that for compact subsets ${\X}\subset\bbc^n$,
$$
\PH {\X} \cap \bbc^n\ =\ \{z\in \bbc^n: \Vaff{\X}(z)\ <\ \infty\}
\tag6.9$$

It is natural to expand $\SHHaffine_{\X}$ by using the  {\bf Lelong class}
$\cl$ of all plurisubharmonic functions $\shaff$ on $\bbc^n$ such that
$\shaff(z)\leq c+\log(1+|z|)$ for some constant $c$.  Set
$$
{\cl}_{\X}\ \equiv\ \{\shaff\in \cl : \shaff  \ \leq\ 0 \  \text{ on}\ {\X}\}
$$   
In analogy with Proposition 4.2 one has  

\Prop{6.2}
$$
\Vaff{\X}(z)\ =\ \sup_{\shaff\in {\cl}_{\X}} \shaff(z)
$$

\vskip .3in

\def\opn{{\co{}}}

\centerline{\bf 7.   Theorems of Sadullaev.}\medskip

The extremal function $\Vaff{\X}$ was studied by A. Sadullaev who proved the
following deep and beautiful result.

\Theorem{7.1 ([S])}  {\sl Let ${\X}\subset Z$ be a compact
subset of an analytic subvariety $Z$ defined in some open subset of 
$\bbc^n$.  Assume that ${\X}$ is not pluripolar in $Z$.  Then if 
$\Vaff{\X}$ (or equivalently ${\V}_{\X}$) is locally bounded on $Z$, the
variety $Z$ must be algebraic. }

\Note{} A subset ${\X}\subset Z$ is {\sl pluripolar} in $Z$  if for each
point $x\in Z$ there is a neighborhood $\opn$ of $x$ in $Z$ and a
plurisubharmonic function $u:\opn\to \bbr$, not identically $-\infty$, such
that  ${\X}\cap \opn\subseteq \{x\in \opn:u(x)=-\infty\}$.

Sadullaev also proved the following.
\Theorem{7.2 ([S])}  {\sl Let $A\subset \bbc^n$ be an irreducible
algebraic curve, and ${\X}\subset A$ a compact subset which is not
pluripolar (equivalently, has positive capacity) in $A$.  Then the extremal
function $\Vaff{\X}$ is harmonic on $A-\X$.

More generally if $A$ is an algebraic subvariety of dimension $m$, then
$\Vaff{\X}$ is the limit of an  increasing sequence of maximal functions
on $A-\X$ (cf. [BT]). }

\medskip
We adapt the arguments of [S] to prove the following   useful result.

\Theorem{7.3}  {\sl  Let $Z$ be a regular  complex analytic curve defined
in some open subset of $\bbc^n$ and consider a compact subset ${\X}\subset
Z$.  Suppose
 $\opn\subset \bbc^n-\X$ is an open set with the property that 
 $$
\Vaff{\X} \equiv \infty \qquad\text{ in }\ \  \opn-Z.
$$
Then in every connected component of $\opn\cap Z$ either 
  $\Vaff{\X}\ \equiv \ \infty$ or $\Vaff{\X}$ is a bounded harmonic function.

In particular if $\Vaff{\X}\equiv \infty$ in $\sim Z$, then in every connected
component of $Z-{\X}$ either $\Vaff{\X}\equiv \infty$ 
or $\Vaff{\X}$ is a bounded harmonic function. }

\pf 
Choose analytic coordinates
$(\zeta, z_1,...,z_{n-1})$ for $|\zeta|<3$ and $\|z\|<1$ on a subset $\opn_0$
of $\opn$ such that $\opn_0\cap Z$ is defined by $z=0$. Let
$D=\{(\zeta,0):|\zeta|\leq 1\}$. 

Recall from Proposition 6.2 that
$\Vaff{\X}(x)= \sup \{u(x): u\in \cl_{\X}\}$. 

\Lemma{7.4} {\sl Fix $u\in \cl_{\X}$ and consider the function $v:Z\to \bbr$
such that $v=u$ in $Z-D$ and $v$ is the harmonic extension of
$u\bigl|_{\partial D}$ to $D$.  Then for every $\epsilon >0$ there exists
a function $u_{\epsilon}\in \cl_{\X}$ whose restriction to $Z$
satisfies             $u_\epsilon \geq \max\{u, v-\epsilon\}$. }

\pf Let $\pi(\zeta,z)=\zeta$ be projection in the coordinate bidisk, and
set  $\wt v = v\circ \pi$. Choose $\rho>0$ sufficiently small that
$$
\wt v -\epsilon\  < \ u\qquad \text{ on the set } \ \{(\zeta,z): |\zeta|=2
\text { and }\|z\|\leq\rho\}.
$$ 
Since $\Vaff{\X}(\zeta,z)=\infty$ for $\zeta\neq 0$, a standard compactness
argument shows that for any $\gamma\in \bbr$ there exists a finite set of
functions $\vf_1,...,\vf_N\in\cl_{\X}$ such that 
$
\vf \equiv \max\{\vf_1,...,\vf_N\} \ >\ \gamma
$ for $ |\zeta|\leq2$ and $\|z\|=\rho $.
In particular we may assume that 
$$
\wt v\ <\  \vf \qquad \text{ on the set } \ \{(\zeta,z): |\zeta|\leq 2
\text { and }\|z\|=\rho\}.
$$
  We now define
$$
u_{\epsilon}\ \equiv\ \cases  \max\{\vf, u, \wt v-\epsilon\}\qquad
\text{for }  \  |\zeta|\leq 2
\text { and }\|z\|\leq\rho\\
\max\{\vf, u\}\qquad\qquad\ \ 
\text{elsewhere } 
\endcases
$$
This proves the lemma.  \qed

To prove Theorem 7.3, suppose first that $\Vaff{\X}$ is bounded on $D$
and let $\{u_m\}_{m=1}^\infty$ be a monotone
increasing sequence from $\cl_{\X}$ converging to $\Vaff{\X}$ on $D$. For each
$u_m$ let $u_{m,\epsilon_m}$  be the function given by Lemma 7.4 with
$\epsilon_m=1/m$. Then on $D$ we have
$$
{u}_m-{\tsize \frac 1m}
\ \leq\ {\wt v}_m-{\tsize \frac 1m} \ \leq \ u_{m,\epsilon_m}
\ \leq\ \Vaff{\X}.
\tag{7.1}
$$
Thus $\Vaff{\X}$ is the limit of the monotone sequence of harmonic functions
$\{{\wt v}_m\}_m$ and must be harmonic.

Suppose now that $\Vaff{\X}$ is unbounded on $D'=\{\zeta\in D:|\zeta|\leq
\frac14\}$.  Then there are sequences $\zeta_m\in D'$ and $u_m\in \cl_{\X}$
with $u_m(\zeta_m)\geq m$.  From (7.1) and the Harnak inequality we
conclude that $\Vaff{\X}\equiv \infty$ on $D'$. This completes the proof.
\qed

\medskip

\Theorem{7.5} {\sl  The last assertion of Theorem 7.3 holds if the
variety  $Z$ has the property that for each point $p\in Z$
there is a neighborhood $\opn$ of $p$ and  an quasi-plurisubharmonic
function $\psi$ on $\bbc^n$ such that  $Z\cap
\opn=\{x\in \opn: \psi(x)=-\infty\}$.    }

\pf
Replace $\vf$ in the argment above with the function $\psi+C$ for
sufficiently large  $C$.

\medskip

\Note{}It is proved in [GZ, Thm. 5.2] that there always exists  such a
$\psi$ with
 $$
Z\cap \opn\subseteq\{x\in \opn: \psi(x)=-\infty\}.
\eqno{(7.2)}
$$
If  equality holds in (7.2) for some $\psi$, $Z$ is called {\sl completely
pluripolar in} $\Omega$. There are   papers which describe how far
a subvariety $Z\subset\Omega$ is from being completely pluripolar. See
[Wie] for example. 

\vfill\eject

\centerline{\bf 8. The Theorem of Fabre.}\medskip
One might hope that projective hulls could be approached by studying 
polynomial hulls in affine open subsets of $\bbp^n$.  This hope is
essentially dashed by the following striking result of B. Fabre, whose
proof uses the generalized Jacobians of singular curves introduced by 
M. Rosenlicht [R].

\Theorem{8.1 (B. Fabre [F$_1$])}  {\sl Let $C\subset \bbp^n$ be an
irreducible algebraic curve with   non-empty singular set. Then there
exist domains $\Omega\subset C$ with smooth boundary having the property
that $\Omega$ meets every divisor in $\bbp^n$.}
\medskip

Note that the resulting Riemann surface with boundary $\Omega\subset
\bbp^n$ has the property that it is not contained in any affine open
subset of $\bbp^n$.

\vskip .3in

\centerline{\bf 9.  Examples}\medskip
There are cases where we understand  the projective hull
completely.

\Prop{9.1}  {\sl  Let $Z\subset \bbp^n$ be an irreducible
algebraic   subvariety of $\bbp^n$, and suppose that ${\X}\subset Z$
is a subset which is not pluripolar (equivalently, has positive capacity)
in $Z$.  Then $\PH {\X} = Z$.}
 
\pf  Choose an affine chart  $Z_0\equiv Z\cap \bbc^n$ and
a compact subset ${\X}_0\subset {\X}\cap Z_0$ which is not pluripolar.
It follows from [S, Prop. 2.1]  that $\Vaff{{\X}_0}$ is locally bounded
on $Z_0$. Therefore, by (6.8), ${\V}_{{\X}_0}$ (and so also ${\V}_{\X}$)
is finite at all points of $Z_0$. Repeating the argument on slight
perturbations of the  chart shows that ${\V}_{\X}$ is finite on all of $Z$, 
i.e., $Z\subseteq \PH {\X}$. However, by Proposition 3.1, we have $\PH
{\X}\subseteq Z$. \qed

\medskip\noindent
{\bf Alternate Proof.}\ \ Let ${\V}_{{\X},Z}$ be the intrinsic extremal
function for  ${\X}$ on the K\"ahler manifold $Z$ (with K\"ahler form
induced from $\bbp^n$).  Guedj and Zeriahi establish Theorem 4.3 for
${\V}_{{\X},Z}$ on $Z$. However, by Proposition 17.4 we have  
${\V}_{\X}\bigl|_Z={\V}_{{\X},Z}$, and the result follows. \qed

\medskip
The projective hull of a subset ${\X}\subset \bbp^n$ is not always an
algebraic set. 

\Theorem {9.2}  {\sl Let   $V=\{(z,f(z))\in \bbc^2 : z\in \bbc\}\subset
 \bbp^2$ be the graph of an entire
holomorphic function $f:\bbc\to \bbc$ which is not a polynomial. Then for
any compact subset ${\X}\subset V$ we have  }
$$
\PH {\X}\ =\ {\X}\cup\{\text{the bounded components of}\ \  V-{\X}\}
$$
\pf
The bounded components of $V-{\X}$ lie in the polynomial hull of ${\X}$
and therefore in the projective hull $\PH {\X}$ by Corollary 3.4.
By Theorem 7.1 the extremal function $\Vaff{\X}$ cannot be locally bounded 
on all of $V$ because $V$ is not algebraic.  We shall prove that
$ 
\Vaff{\X}\equiv-\infty
$ 
in $\bbc^2- V$. Then by Theorem 7.3 we must have $\Vaff{\X}\equiv
-\infty$ on the unbounded component of $V-{\X}$ and the assertion will
follow.

It remains to prove that $\PH {\X}\subset V$.  By applying a
homothety to $\bbc^2$ we can assume that  $\pi({\X})\subset \{z\in \bbc:
|z|\leq 1/2\}$ where $\pi:\bbc^2\to\bbc$ denotes   projection onto the
first coordinate.  Since $f$ is entire, it has power series expansion
$f(z)=\sum_{n=0}^\infty a_nz^n$ with $$
{\limsup_{n\to \infty}} \,|a_n|^{\frac1n}\ =\ 0.
\tag9.1$$
Fix $(z_0,w_0)\notin V$ and consider now the family of polynomials
$
P_d(z,w) = w-\sum_{n=0}^d a_nz^n
$
for $d\geq 1$. Since $w_0\neq f(z_0)$ we have 
$$
|P_d(z_0,w_0)|\geq  \frac12|w_0-f(z_0)|>0   \qquad\qquad
\text{for all $d$ sufficiently large.} 
\tag9.2
$$
 However, since ${\X}$ is contained in the graph of $f$ over
$\{|z|\leq\frac12\}$, we have 
$$
\sup_{\X}|P_d|\leq \sup_{|z|\leq\frac12}\left\{\sum_{n=d+1}^\infty
|a_n||z|^n\right\} <  \sup_{n>d} \,|a_n| = 
\left\{\sup_{n>d}\,|a_n|^{\frac1d}\right\}^d \leq 
\left\{\sup_{n>d}\,|a_n|^{\frac1n}\right\}^d
 $$
for large $d$, and so for any $C>0$ we have 
$$
\lim_{d\to \infty} C^d\sup_{\X}|P_d|\leq \lim_{d\to \infty}C^d
\left\{\sup_{n>d}\,|a_n|^{\frac1n}\right\}^d=
\lim_{d\to \infty}\left\{C \sup_{n>d}\,|a_n|^{\frac1n}\right\}^d
\leq C \limsup_{n\to \infty}\,|a_n|^{\frac1n}.
$$
By (9.1) the rightmost term is zero. Hence there exists no constant
$C>0$ such that $|P_d(z_0,w_0)|\leq  C^d\sup_{\X}|P_d|$ for all $d$, and so 
by Proposition 6.1    $(z_0,w_0)\notin \PH {\X}$.
 \qed

\medskip

Further interesting examples come from classical gap series.
We recall the following.

\Theorem{9.3 (See Hille [Hi])} {\sl Consider the holomorphic function
$$
f(z)=\sum_{k=0}^\infty c_{n_k}z^{n_k}
\qquad\text{ with }\qquad 
\limsup_{k\to \infty}\left\{|c_{n_k}|^\frac{1}{n_k}\right\}=1.
\tag9.3
$$
Assume that there exists $\lambda >1$ with $\lambda n_k < n_{k+1}$
for all $k>k_0$.  Then $\Delta_1(0)=\{|z|< 1\}$ is the domain of
analyticity of $f$ }

\Theorem{9.4} {\sl Let $f$ be as in Theorem 9.3 and assume   the series
(9.3) converges for all $|z|=1$.  Fix any $r$, $0<r<1$, and let 
$$
V_r\ \equiv \ \{(z,f(z)\in \bbc^2 : |z|< r\}
\and
dV_r\ \equiv \ \{(z,f(z)\in \bbc^2 : |z|= r\}
$$
If
$$
\limsup_{k\to\infty}\frac{n_{k+1}}{n_k} \ =\ \infty,
\tag9.4
$$
then for all $r$, 
$$
\PH{dV_r}\ =\ V_r
$$
}

\pf As in the proof of Theorem 9.2 it will suffice to prove that
$\PH{dV_r}\subseteq V_1$. To see this,
fix  $z_0\in \overline{\Delta}_1$  and
choose $w_0\neq f(z_0)$. Set $
P_{n_k}(z,w) = w-\sum_{j=1}^k c_{n_j} z^{n_j}
$
and write 
$
P_{n_k}(z,w) = w_0-f(z_0)+ \sum_{j=k+1}^\infty c_{n_j} z_0^{n_j}.
$
Then for all $k$ sufficiently large we have
$$
|P_{n_k}(z,w)|\ \geq\ |w_0-f(z_0)| - 
\left| \sum_{j=k+1}^\infty c_{n_j} z_0^{n_j} \right|
\ \geq\  
 >\ 0.
$$
On the other hand, using (9.3) above we find that
$$\aligned
\sup_{dV_r}\ &=\ 
\sup_{|z|=r} \left|f(z)- \sum_{j=1}^k c_{n_j} z^{n_j}  \right| \\
&=\ \sup_{|z|=r} \left| \sum_{j=k+1}^\infty c_{n_j} z^{n_j}  \right| 
\leq \ \sum_{j=k+1}^\infty\left(|c_{n_j}|^{\frac1{n_j}} r \right)^{n_j} \\
&\leq \ \sum_{j=k+1}^\infty\left(\sup_{\ell\geq j}
|c_{n_\ell}|^{\frac1{n_\ell}} r  \right)^{n_j} 
\leq \ \sum_{j=k+1}^\infty (\sqrt r )^{n_j} \ <\ L r^{\frac{n_{k+1}}2}
\endaligned
$$
for $L=(1-\sqrt r)^{-1}$ and for all $k$ sufficiently large.

Now if $(z_0,w_0)\in \PH {dV_r}$, then by Proposition 6.1 and the
paragraph above, there exists $C>1$ with

$$
K\ \leq\ C^{n_k} L r^{\frac{n_{k+1}}2}.
$$
Hence, there exists $K_0>0$, $a>0$ and $b>0$ with
$$
K_0 \ \leq \ \exp({an_k-bn_{k+1}})
\qquad\text{ for all $k$ sufficiently large.}
$$
Let $\kappa_0=\log(K_0)$.  Then 
$$
\kappa_0 \ \leq \  an_k-bn_{k+1}
\qquad\text{ for all $k$ sufficiently large.}
$$
or equivalently
$$
\frac{n_{k+1}}{n_k} + \frac{\kappa_0}{b n_k}\ \leq\ \frac a b
$$
which contradicts assumption (9.4).  We conclude that $(z_0,w_0)\notin
\PH{dV_r}$.

Suppose now that $|z_0|=R>1$ and $w_0$ is arbitrary.
 Set $
P_{n_k}(z,w) = w-\sum_{j=1}^k c_{n_j} z^{n_j}
$
as before and note that by (9.3)
$$\aligned
|P_{n_k}(z_0,w_0)|\ &\geq\ |c_{n_k}|R^{n_k} -
|w_0|-\sum_{j=1}^{k-1}|c_{n_j}z_0^{n_j}|  \\
  &\geq\ \left(|c_{n_k}|^{\frac{1}{n_k}}R \right)^{n_k}-
|w_0|- K\sum_{j=1}^{k-1}\left(\rho R\right)^{n_j}
\endaligned
$$
for all $k$ is sufficiently large,   where $\rho$ satisfies
$1<\rho< R$. Choose   $\a<1$ with $\rho <\a^2 R$.
Then for $k$ sufficiently large
$$\aligned
   \left(|c_{n_k}|^{\frac{1}{n_k}}R \right)^{n_k}-
|w_0|- K\sum_{j=1}^{k-1}\left(\rho R\right)^{n_j} 
&\geq\  \left(\a  R \right)^{n_k}- K_1\sum_{m=0}^{n_{k-1}}(\rho R)^m\\
&=\ \left(\a  R \right)^{n_k}- K_1\frac{(\rho R)^{n_{k-1}+1}-1}{\rho
R-1}          \\ &\geq\  \frac{(\a  R)^{n_k}(\rho R-1)-K_1(\rho
R)^{n_{k-1}+1}}{\rho R-1}. \endaligned
\tag9.5
$$
The numerator of the last term in (9.5) is
$$\aligned
\left[\left(\frac{\a}{\rho}\right)^{n_{k-1}+1}(\a R)^{n_k-n_{k-1}-1}   
 (\rho R-1)-K_1\right]&(\rho R)^{n_{k-1}+1}   \\ 
&\geq\ \left(\frac{\a}{\rho}\right)^{n_{k-1}}(\a R)^{n_k-n_{k-1}-1} 
(\rho R-1)-K_1\\
&\geq\ \left(\frac{\a^2R}{\rho}\right)^{n_{k-1}}(\a R)^{n_k-2n_{k-1}-1}
(\rho R-1)-K_1\\
&\geq\ (\a R)^{n_k-2n_{k-1}-1}
(\rho R-1)-K_1\ \  \arr\ \  \infty
 \endaligned
$$ 
since $\a R>\a^2 R/\rho>1$ and $n_k-2n_{k-1}\ \to\ \infty$.
In particular $|P_{n_k}(z_0,w_0)|\geq K_2>0$ for all $k$ sufficiently
large, and  our previous estimate for $\sup_{dV_r}|p_{n_k}|$ rules 
this case out as well.\qed

\Remark{9.4. (A Projectively Convex Curve)} We claim that for the curve
$$
\G\ =\  \left\{(z,\exp\left(z+\overline{z}\right) :
 |z|=1\right\} \ \subset\ \bbc^2\ \subset\
\bbp^1\times \bbp^1\subset \bbp^3, \qquad \text{one has }\ \  \PH {\G}=\G 
$$
An outline of the proof is as follows. Arguing as in Theorem 9.2 one sees
that $\PH{\G}$ is contained in 
$$
Z\ =\  \left\{(z,w) : w=\exp\left(z+\frac 1z\right), 
\ \ 0<|z|<\infty\right\}
\cup (0\times \bbp^1) \cup (\infty\times\bbp^1 ).
$$
By Theorem 7.3 the extremal function $\Lambda^0_{\G}$ (and therefore also
$\Lambda_{\G}$) is either $\equiv \infty$ or is a locally bounded
function on each of the two  components of $Z$ over $\{0<|z|<1\}$ and  
$\{1<|z|<\infty\}$.  However, the automorphism $(z,w)\mapsto(1/z,w)$
shows that either  $\Lambda^0_{\G}\equiv\infty$ on both components or
it is bounded on both. However, by Sadullaev' Theorem 7.1 it cannot be
bounded on both, and therefore
$$
\PH \G\ \subseteq\ \G
\cup (0\times \bbp^1) \cup (\infty\times\bbp^1 ).
$$
Suppose now that there exists a point $x\in {\PH \G}\cap (0\times
\bbp^1)$. By Theorem 11.1 there exists a probability measure $\mu$ on
$\G$ and a positive current $T$ of bidimension (1,1) with support in 
$\PH \G$, such that $dd^c T=\mu-\delta_x$. This is impossible since, if
$\pi:\bbp^1\times\bbp^1\to\bbp^1$ denotes projection on the first factor,
we would have $dd^c \pi_* T = \pi_*\mu-\delta_0$ with $\supp(\pi_* T)$ in
the unit circle $S^1$ and $\pi_*\mu$ a probablilty measure on $S^1$.
By symmetry we conclude that $\PH \G = \G$.

\vskip .3in

\centerline{\bf 10. Compactness and Stability.}\medskip

While we have many representations of the projective hull, we do not
  have an easy answer to the following natural question: {\sl  Given a
compact subset ${\X}\subset \bbp^n$, is $\PH {\X}$ also compact}?
Proposition 9.1  shows that if ${\X}$ is contained in an affine chart
$\bbc^n$, its projective hull may not be contained in that chart.
For example consider the tiny curve 
$$
{\X}=\{(t, (1+t)^\frac1m)\in\bbc^2 : |t|=\epsilon\}
$$
for some choice of $m$th root and $\epsilon$ very small. 
Then $\PH {\X} \cap \bbc^2=\{(x,y):y^m=1+x\}$   is an algebraic curve
which is $m$-sheeted over the entire $x$-axis.

\Def{10.1}  A subset ${\X}\subset \bbp^n$ is called {\sl stable} if the
best constant function $C:\PH {\X}\to \bbr^+$ is bounded above, or
equivalently, the radius function $\rho_{\X}:\PH {\X}\to \bbr^+$ defined in
(5.4) is bounded below by a positive constant (cf. Proposition 5.2). 
\medskip

Note that $\PHull {\SH {\X}}$ has an {\sl outer boundary} consisting of all
points $Z\in \PHull {\SH {\X}}-\{0\}$ with $\|Z\|=\rho(\pi Z)$. Stability
is equivalent to the fact that $0$ is not in the closure of this outer
boundary.

\Prop{10.2 }\ \ \  {\sl If ${\X}\subset \bbp^n$ is stable, then $\PH {\X}$
is compact.}

\pf  
If $C_{\X}(x)\leq C<\infty$ for all $x\in\PH {\X}$, then for every section
$\cp$ of $\co {}(d)$
$$
\|\cp(x)\|\ \leq\ C^d\sup_{\X}\|\cp\|\qquad\text{ for all } x\in\PH {\X}
$$
and hence for all $x$ in the closure of $\PH {\X}$, proving that $\PH
{\X}$ is closed.\qed

\medskip

We next examine some elementary local properties of
projective hulls.   

Fix a closed set ${\X}\subset \bbp^n$.

\Prop{10.3}  {\sl  Suppose $B\subset \bbp^n$ is a closed subset with
$C_{\X}$   bounded above on $\PH {\X} \cap B$. Then $\PH {\X}\cap B$ is
compact and  $$ (\PH {\X}\cap B)\widehat{\ \ }\cap B\  =\ \PH {\X} \cap B
$$
In particular if $\PH {\X}$ is stable, then}
$$
\PH {\PH {\X}}\   =\ \PH {\X} 
$$

\pf  The compactness is proved as in 10.2. For the second assertion 
it suffices to show that
 $(\PH {\X}\cap B)\widehat{\ \ }\ \subseteq\ \PH {\X}$. By assumption 
 there exists a constant $C<\infty$ with
 $$
\sup_{\PH {\X} \cap B}\|\cp\| \ \leq C^d \sup_{\X} \|\cp\|
$$
for all $\cp\in H^0(\bbp^n,\co{}(d))$ and all $d\geq 0$. If 
$x\in (\PH {\X}\cap B)\widehat{\ \ }$, there exists $C_x$ such that
$$ 
\|\cp(x)\|\ \leq \ C_x^d \sup_{\PH {\X}\cap B} \|\cp\|   
 \  \leq(C_x C)^d \sup_{\X} \|\cp\|
$$
for all $\cp, d$ as above. Hence, $x\in \PH {\X}$.\qed
\medskip

We now switch to affine coordinates.

\Prop{10.4} {\sl Suppose $B$ is a compact polynomially convex subset
of $\bbc^n\subset \bbp^n$.  If $C_{\X}$ is bounded on $\PH {\X}\cap B$,
then  $\PH {\X} \cap B$ is  polynomially convex.}

\pf Set 
$$
C_B\ \equiv \ \sup_{z\in \PH {\X}\cap B} C(z)(1+|z|^2)^{\frac12}\
<\ \infty.
$$
Then it follows directly  from (6.3) that  for any polynomial
$p\in\bbc[z_1,...,z_n]$ of degree $d$
$$
\sup_{\PH {\X}\cap B}|p(z)|\ \leq C_B^{d} \sup_{\PH {\X}}\|\cp\|
$$
where $\cp$ is the associated section of $\co{}(d)$.
Now if $x\in \PHull{[\PH {\X} \cap B]}$, then
$$
\frac{|p(z)|}{(1+|z|^2)^{\frac d2}}\ \leq\ 
|p(z)|\ \leq\ \sup_{\PH {\X} \cap B} |p| \ \leq\ C_B^d  \sup_{\PH {\X}}\|\cp\|
$$
for all polynomials $p$ as above, and so $x\in \PH {\X}$. Furthermore,
$$
|p(z)|\ \leq\ \sup_{\PH {\X} \cap B} |p| \ \leq\ \sup_B |p|
$$
and so $x\in \PHull B= B$.  Hence $x\in \PH {\X}\cap B$.\qed

\vfill\eject

\centerline{\bf 11. Projective Hulls and Jensen Measures -- The
Main Result}\medskip

 In this section we shall prove our central analytic result concerning
projective hulls. A general version of the theorem will be established in
\S 18, but here we   shall work on $\bbp^n$.  Let 
$$
P\ \equiv\ \PSH_{\omega} = \{\vf\in C(\bbp^n) : dd^c \vf
 +\omega  \geq 0 \} 
$$
be the space of continuous quasi-plurisubharmonic functions where
$\omega $ is the standard K\"ahler form. Then the projective hull
of a compact subset ${\X}\subset \bbp^n$ is the set $\PH {\X}$ of points
$x$ for which there exists a constant $\lambda=\lambda_x$ such that  
$$
\vf(x) \leq \lambda+ \sup_{\X} \vf \qquad\text{ for all } \vf \in P
\tag11.1$$
Denote by $\PH {\X}(\Lambda)$ the subset of $\PH {\X}$ for which 
there exists a $\lambda\leq \Lambda$ satisfying condition (11.1).
Let
$\cp_{1,1}(\Lambda)$ denote the convex cone of positive currents  of
bidimension $(1,1)$ and mass $\leq \Lambda$ on $\bbp^n$, and $\cm_{\X}$
the set of probability measures on $\bbp^n$ with support in ${\X}$.

\Theorem {11.1}  {\sl For a compact subset ${\X}\subset \bbp^n$ the
following are equivalent:

\medskip

\qquad (A)\ \ \ \  $x\in \PH {\X}(\Lambda)$
\medskip

\qquad (B) \ \ There exist  $T\in \cp_{1,1}(\Lambda)$  and
$\mu\in\cm_{\X}$ such that:\medskip 

 \hskip 1.6in (i)    \ \ \ $dd^c T\ =\ \mu - \d_x$
\smallskip

 \hskip 1.6in (ii)    \ \ $\supp(T)\ \subset {\PH {\X}}^{-}$ \smallskip
}

\ 
\pf 
What follows is a succinct proof of the result.  In \S 18 a more rounded
and geometric proof is given of a more general result. The
interested reader may want  to go directly there.

We first show that (B) $\Rightarrow$ (A).  Suppose
$T\in \cp_{1,1}(\Lambda)$ satisfies (i) and (ii).
  Let $\vf\in \PSH_{\omega}$, so that 
$$
dd^c \vf\ + \omega \ =\ \eta\ \geq\ 0
$$
 Then
$$\aligned
\int_{\X}\,\vf\, d\mu - \vf(x)\ &=\ (dd^cT)(\vf)\ =\ T(dd^c\vf) 
 =\ T(\eta)-T(\omega)   \\
&\geq\ -T(\omega)\ =\ - M(T)\ \geq\ -\Lambda
\endaligned
$$
Hence, 
$
\vf(x)\ \leq \ \sup_{\X} \vf +\Lambda
$, and we conclude that $x\in \PH {\X}(\Lambda)$.
\medskip
To show that (A) $\Rightarrow$ (B) we shall need the following.

\Lemma {11.2}  {\sl Fix a compact subset $B\subset \bbp^n-{\PH {\X}}^-$
and real numbers $M,N$ with $N-M>0$.  Then for all $\Lambda>0$ there exists 
a function $\vf\in P$ such that 
$\vf\bigl|_{\PH {\X} (\Lambda)}\leq M$ and $\vf\bigl|_B\geq N$.
}

\pf Condition (1) and compactness enables us to find a finite 
set $\vf_1, ..., \vf_m\in P$ such that $\sup_{\X}\vf_i\leq M-\Lambda$  for all
$i$ and  $\vf\equiv\max\{\vf_1,...,\vf_m\}\geq N$ on $B$.  \qed
\medskip

We now fix $\epsilon >0$ and consider the
closed convex cone 
$$
\cp_{1,1}(\Lambda, \epsilon) \equiv 
\{ T\in \cp_{1,1}(\Lambda) : \supp(T)\subset {\PH {\X}}_{2\epsilon}\}.
$$
where  
$$
{\PH {\X}}_t = 
\{x\in \bbp^n : \text{dist}(x, \PH {\X})\leq t\}.
$$

It will suffice to show  that there exists a current $T\in
\cp_{1,1}(\Lambda, \epsilon) $ satisfying conditions (i) and (ii) above.
If this is not the case, i.e., if 
$$
\left(\cm_{\X} - \delta_x\right) \,\cap\, \cp_{1,1}(\Lambda, \epsilon) 
\ =\ \emptyset
$$
then by the Hahn-Banach Theorem there exists $\vf \in C^{\infty}(\bbp^n)$
and $\gamma\in \bbr$ such that
$$
\sup_{\mu\in \cm_{\X}}(\mu-\delta_x)(\vf) \  < \ \gamma \ <  (dd^cT)(\vf)
$$
for all  $T\in \cp_{1,1}(\Lambda, \epsilon) $.
Note that since $0\in \cp_{1,1}(\Lambda, \epsilon)$ we have $\gamma <0$.
Setting $\psi =\frac {\Lambda}{|\g|}\vf$, we find that
$$
\sup_{\mu\in \cm_{\X}}(\mu-\delta_x)(\psi) \ 
<\ -\Lambda\ < (dd^cT)(\psi) = T(dd^c\psi)
\tag11.2$$
for all  $T\in \cp_{1,1}(\Lambda, \epsilon)$.
Applying the right hand inequality to currents of the form
$T=\delta_y \xi$ where $y\in {\PH {\X}}_{2\epsilon}$ and $\xi$ is a positive
simple (1,1) vector of length $\Lambda$ at $y$, we conclude that
$$
dd^c \psi +\omega \ \geq\ 0 \qquad\text{ on } \ {\PH {\X}}_{2\epsilon}
$$
Now let $\vf_{M,N}\in P$ be the function given by Lemma 2 with $M$ and $N$
chosen so that
$$
M
 <\ \inf_{\PH {\X}(\Lambda)} \psi
\and
N \ >\ \sup_{\bbp^n-{\PH {\X}}_{\epsilon}^0} \psi.
$$
Then the function
$$
\Psi \ \equiv\ \max\{\psi, \vf_{M,N}\}
$$
has the property that 
$$
\Psi=\psi \ \ \text{ on } \PH {\X}(\Lambda)
\and
\Psi = \vf_{M,N} \ \ \text{ on } \bbp^n-{\PH {\X}}_{\epsilon}^0
\tag11.3$$
Consequently, $dd^c\Psi + \omega \geq 0$ on all of $\bbp^n$, that is,
$\Psi \in P$.
Furthermore, by (11.2) and (11.3) we have that
$$
\sup_{\mu\in \cm_{\X}}\int_{\X}\Psi d\mu - \Psi(x) \ < \ -\Lambda
$$
for all probablility measures $\mu$ on ${\X}$. Choosing $\mu=\delta_y$
for $y\in {\X}$ shows that
$$
\sup_{\X}\Psi - \Psi(x) \ <\ -\Lambda
$$
which means that $x\notin \PH {\X}(\Lambda)$.\qed

\Cor{11.3} {\sl  For any   $\nu\in \cm_{\PH {\X} (\Lambda)}$ 
there exists $T\in\cp_{1,1}(\Lambda)$ and $\mu\in \cm_{\X}$ with
\medskip

 \hskip 1.6in (i)    \ \ \ $dd^c T\ =\ \mu - \nu$
\smallskip

 \hskip 1.6in (ii)    \ \ $\supp(T)\ \subset {\PH {\X}}^{-}$ \smallskip
}

\pf The probability measures on $\PH {\X}  (\Lambda)$ are the closed
convex hull of the $\delta$-measures.\qed

\Cor{11.4} {\sl  For any  $x\in {\PH {\X}}^-$  there  are probability
measures $\nu\in \cm_{{\PH {\X}}^-}$, $\mu\in \cm_{\X}$ and a current 
$T\in\cp_{1,1}$  
 with \medskip

 \hskip 1.6in (i)    \ \ \ $dd^c T\ =\ \mu - \nu$
\smallskip

 \hskip 1.6in (ii)    \ \ $x\in \supp(T)\ \subset {\PH {\X}}^{-}$ \smallskip
}

\pf
Let $\{x_k\}_{k=1}^\infty\subset \PH {\X}$ be  a sequence converging to $x$.
Choose currents $T_k\in \cp_{1,1}$ as in Theorem 11.1 with 
$dd^cT_k=\mu_k-\delta_{x_k}$ and $\supp(T_k)\subset {\PH {\X}}^{-}$.
Then the positive current
$$
\wt T\ \equiv \ \sum_{k=1}^\infty \frac1{2^kM(T_k)} T_k
$$
has $x\in \supp(T)$ and satisfies $dd^c\wt T=\wt \mu-\wt \nu$ for positive
measures $\wt \mu $ on ${\X}$ and $\wt\nu$ on ${\PH {\X}}^-$.
$T\equiv\frac 1{\mu({\X})}\wt T$ is the desired current.\qed 

\medskip
\Remark{11.5} Let $\Lambda_K(x)$ be the extremal function
introduced in (1.2) and discussed in \S 4.  Note that by definition:
$$
x\in \PH K (\Lambda)\ \ \Leftrightarrow \ \ \Lambda_K(x)\leq \Lambda
$$
For a fixed point $x\in \PH K$ let $\cf_x$ denote the set of positive
currents $T$ of bidimension (1,1) satisfying $dd^cT=\mu-\delta_x$ for
some $\mu\in \cm_K$ and set ${\cn_K}=\{T\in \cm_K: \supp(T)\subseteq       
{ \PH K}^-\}$. Let $T_x$ be the current guaranteed by Theorem 11.1.
Then from the discussion above we have that
$$
\Lambda_K(x)\ =\ M(T_x)\ =\ \inf_{ \cm_K} M(T) \ =\ \inf_{ \cn_K} M(T). 
\tag11.4
$$
In other words {\sl $T_x$ is the positive  current of least mass satisfying
the equation $dd^cT=\mu-\delta_x$, and that least mass is exactly
$\Lambda_K(x)$.}

The middle equality in (11.4) follows from the fact that Theorem 11.1
also holds without the requirement  $\supp(T)\subseteq  { \PH K}^-$
as one can easily check (cf. Theorem 18.2 below).

\Remark{11.6} Suppose that $K=\gamma$ is a closed curve and $\PH K =V$
is a 1-dimensional complex submanifold with boundary
$\gamma$ as in the examples below.  Then
$$
T_x\ =\ G_x [V]
$$
where $G_x$ is the Green's function on $V$ with singularity at $x$.

\vskip.3in 
\vfill\eject
\centerline{\bf 12. Structure Theorems}

\medskip

Theorem 11.1 has a number of basic consequences.  We recall that a subset
$W$ of a complex manifold $Z$ is called {\bf 1-concave} if for every open
set $\co{}\subset\subset Z$ and every holomorphic map $f$ from a
neighborhood of $\overline{ \co{}}$ to $\bbc$ one has 
$f(W\cap\co{})\subset\bbc-\Omega$ where $\Omega$ is the unbounded
component of $\bbc-f(W\cap \partial \co{})$.

Fix a compact subset ${\X}\subset \bbp^n$

\Theorem{12.1} {\sl The set    
${\PH {\X}}^- - {\X}$ is 1-concave in $\bbp^n-{\X}$.}

\Note{}    {\sl It follows that
$({\PH {\X}}^- - {\X})\cap O$ is 1-concave in $O-{\X}$
 for any open subset} $O\subset \bbp^n$.

\pf
Corollary 11.4 implies that for any $x\in {\PH {\X}}^- -{\X}$ there exists
a positive $(1,1)$-current $T$ with 
$$\aligned
&(i)\ \ \ \ dd^cT = -\nu \ \leq\ 0\ \ \ \text{ in } \bbp^n-{\X}
\qquad\text{ and } \\
&(ii)\ \ \ \ x\in \supp(T)\subseteq {\PH {\X}}^-.
\endaligned
\tag12.1
$$
By  Proposition 2.2 of [DL] we conclude that
supp$(T)$ is 1-concave in $\bbp^n-{\X}$. 
 The argument given for Proposition 2.3 of  [DL] now gives the following:
\Lemma {12.2}. {\sl Let $\cs$ denote the family of (closed) 1-concave
subsets of $\bbp^n-{\X}$ which have support in ${\PH {\X}}^-$.   Then the
union of all elements in $\cs$ is again an element of $\cs$.}

\medskip\noindent
By (12.1)  this maximal 1-concave subset  
 equals  ${\PH {\X}}^-$. \qed

\medskip

\Cor{12.3} {\sl For any compact subset ${\X}\subset \bbp^n$ the
closed projective hull ${\PH {\X}}^- -{\X}$ has locally positive Hausdorff
2-measure.}

\pf Any 1-concave subset of a complex manifold has positive
Hausdorff 2-measure in any neighborhood of any point.
(Otherwise the complement of the image under any holomorphic map to $\bbc$
would be connected.) \qed

\Theorem {12.4}  {\sl Let ${\X}\subset \bbp^n$ be any compact subset and
assume that the Hausdorff 2-measure of ${\PH {\X}}^-$ is locally finite. 
Then ${\PH {\X}}^-  -{\X}$ is a 1-dimensional analytic subvariety of
$\bbp^n-{\X}$.

Moreover, suppose $D\subset\bbc^n-{\X}\subset\bbp^n$ is a strictly convex
domain with smooth boundary. If the Hausdorff 1-measure of ${\PH {\X}}^- 
\cap \partial D$ is finite, then ${\PH {\X}}^- \cap D$ is a 1-dimensional
analytic subvariety of $D$.
}

\pf  This follows from Theorem 12.1 and [DL, Thm. 3.3 and Cor. 3.8].
It also follows from Theorem 5.7 and [Sib, Thm. 17].\qed

\medskip
\Theorem {12.5} {\sl Let $\Gamma\subset \bbp^2$ be a finite union of 
smooth closed
 curves which are pluripolar. Then the local Hausdorff
dimension of $\PH {\Gamma}$ is  everywhere 2.  The conclusion also holds for any
$\Gamma\subset \bbp^n$ which is real analytic. }

\pf  Since $\Gamma$ is  pluripolar,  $\PH {\Gamma}$ is pluripolar
by Corollary 4.4.
Suppose   $\ch^{2+\a}(\PH \Gamma)>0$ for some
$\a>0$ where $\ch^\b$ denotes Hausdorff measure in  dimension $\b$.
Applying the coarea formula [F,3.2] to a linear projection 
$\pi: \bbp^2-\bbp^0\to \bbp^1$ shows that the set of $y\in \bbp^1$ with 
$\ch^{\a}(\PH \Gamma\cap\pi^{-1}(y))  >0$ has positive $\ch^2$-measure.

Now the condition  $\ch^{\a}(\PH \Gamma\cap\pi^{-1}(y))  >0$ implies that
$\pi^{-1}(y)\subset \PH \Gamma$.   To see this consider the sets $\PH
{\G}_{y,t}\equiv \{x\in  \PH \Gamma\cap\pi^{-1}(y) : C_{\G}\leq
t\}$, and note that $\ch^{\a}(\PH {\G}_{y,t})  >0$ for  $t$ sufficiently 
large. For such $t$,  ProjHull$(\PH {\G}_{y,t})= \pi^{-1}(y)$ since sets of
capacity zero  in $\bbp^1$ have measure zero. However, ProjHull$(\PH
{\G}_{y,t}) \subset \PH \G$, and so $\pi^{-1}(E)\subset \PH \G$ for some
set $E$ of positive measure in $\bbp^1$. Hence, $\ch^4(\PH \G)>0$ and
therefore  $\PH \G = \bbp^2$ contradicting the pluripolarity of $\PH \G$.

To prove the second statement note that any real analytic curve is
pluripolar and that under projections $\pi:\bbp^n-\bbp^{n-3}\to\bbp^2$
one has $\pi(\PH \G)\subseteq \PH{\pi(\G)}$
\qed

\medskip

\Theorem {12.6}{\sl Suppose that 
${\X}\subset \subset\bbc^n$ and let $\PH {\X}_0$
be a connected component of ${\PH {\X}}^- -{\X}$  which is bounded in
$\bbc^n$.  Then $\PH {\X}_0$ is contained in the polynomial hull of ${\X}$.}

\pf
Since ${\PH {\X}}^- -{\X}$ is a 1-concave subset of $\bbp^n-{\X}$, so is
any connected component $\PH {\X}_0$. We now use the fact ([DL, Prop.2.5])
that the polynomial hull of ${\X}$ in $\bbc^n$ is the union of all bounded,
1-concave subsets of $\bbc^n-{\X}$.\qed

\Cor{12.7} {\sl
If $\PH {\X} \subset \Omega=\bbp^n-D$ for some algebraic hypersurface $D$,
then $\PH {\X} ={\PH {\X}}_{\Omega}$.  In particular if $\g\subset\Omega$
is a C$^1$-curve and if $\PH \g\subset\subset \Omega$, then $\PH \g -\g$ 
is a 1-dimensional analytic subvariety of $\Omega-\g$.
}
\pf By Proposition 3.2   taking projective hulls commutes with Veronese
embeddings. However, by embedding $\bbp^n\subset \bbp^N$ by the $d$th
Veronese map, where  $d=$ deg$(D)$, we reduce to the situation of Theorem
12.6.  \qed

\medskip

The next result is a   strong form of the Local Maximum Principle
for projective hulls.

\Theorem {12.8} {\sl Fix } $U^{\text{open}}\subset\subset
\Omega^{\text{affine open}}\subset \bbp^n$.  {\sl Then}
$$\boxed{
{\PH {\X}}^- \cap U \ \subseteq \text{ PolyHull}
\left\{({\PH {\X}}^-\cap\partial
U)\cup ({\X}\cap U)\right\}
}
 $$
{\sl with equality if $C_{\X}$ is bounded on $\PH {\X} \cap \overline U$.}

\pf Set 
$$
\Sigma \ \equiv\ {\PH {\X}}^- - {\X}
$$
$$
O \ \equiv\ \bbp^n - ({\PH {\X}}^- \cap \partial U)
$$
By the Note following Theorem 12.1, $\Sigma\cap O$ is 1-concave in $O-{\X}$.

\medskip
\noindent
{\bf Claim:}\ \ $\Sigma\cap U$ is a union of connected components of 
$\Sigma\cap O$.

\pf  Set $\Sigma_O\equiv \Sigma\cap O$.  Note that
$$
\bbp^n \ =\ U\ \amalg\ \partial U\ \amalg\ (\sim \overline U)
$$
gives
$$\aligned
\Sigma_O \ &=\ (\Sigma_O\cap U)\ \amalg\ 
(\Sigma_O\cap \partial U)\ \amalg\ (\Sigma_O\cap \sim \overline U)\\
&=\ (\Sigma_O\cap U)\ \amalg\ 
 (\Sigma_O\cap \sim \overline U)
\endaligned
\tag12.2$$
since
$$\aligned
\Sigma_O \cap \partial U \ &=\ (\Sigma \cap O)\cap \partial U  \\
&=\ ({\PH {\X}}^- - {\X} {\X})\cap (\bbp^n-({\PH {\X}}^- \cap\partial
U))\cap\partial U \\ 
&=\ \{{\PH {\X}}^- - {\X}-({\PH {\X}}^- \cap \partial U)\}\cap
\partial U \\ &=\ \emptyset.\endaligned
$$
Therefore (12.2) gives a disconnection of $\Sigma_O$. This proves the
claim. \medskip

Now by [DL] we know that for any compact subset $C\subset \Omega$ we have
that the polynomial hull of $C$ equals the 1-Hull of $C$ which is by
definition the union of all bounded  1-concave subsets of $\Omega-C$.

For the last statement recall from Proposition 10.3  that if $C_{\X}$ is
bounded on $\PH {\X}\cap \overline  U$, then  
$(\PH {\X}\cap \overline U )\widehat{\ \ }\cap \overline U \  
=\ \PH {\X} \cap \overline U$.  Now use the fact that for bounded subsets
of $\Omega$, the polynomial hull is contained in the projective hull.
\qed

\medskip
This enables us to give the following  generalization of Wermer's Theorem.

\Theorem{12.9} {\sl Let $\gamma\subset\bbp^n$ be a finite union of real
analytic curves. Then $\PH{\gamma}$ is a subset of Hausdorff
dimension 2 whose closure is  1-concave. Furthermore,  if the Hausdorff
2-measure of ${\PH{\gamma}}^-$ is finite in a neighborhood of some complex
hypersurface, then ${\PH{\g}} -\g ={\PH{\g}}^- -\g$ is a 1-dimensional
complex analytic subvariety of $\bbp^n-\g$.

The same conlcusion holds for any smooth pluripolar curve $\g$  in
$\bbp^2$.}

\pf Theorems 12.1 and 12.5 give the first statement. Suppose now that
$\ch^2({\PH{\g}}^- \cap \co{})<\infty$    where $\co{}$ is a neighborhood
of some divisor $D$. We may assume that $D\cap\g=\emptyset$ and therefore
that $\co{}\cap\g=\emptyset$.  By Theorem 12.4 we know that ${\PH{\g}}^-
\cap \co{}$ is a 1-dimensional complex analytic subvariety of $\co{}$. We
now choose a bounded subdomain $U\subset\subset \Omega = \bbp^n-D$ with
real analytic boundary $\partial U\subset \subset \co{}$. 
Then $\Gamma \equiv({\PH{\g}}^-  \cap\partial U)\cup\g$ is a real analytic
curve (which we may assume to be regular by appropriate choice of $U$),
and by Wermer's Theorem $\PH \Gamma$ is a bounded 1-dimensional complex
subvariety of $\Omega- \Gamma$ with regularity at the boundary as in
[HL$_1$]. In particular it is regularly and analytically immersed up to
the boundary in $\co{}$.   We conclude that  $W\equiv\PH{\G}\cup
(\PH{\g}^-\cap\co{})-\g$ is a 1-dimensional complex subvariety in 
$\bbp^n-\g$.  

 Now by Theorem 12.8 we have $\PH{\g}^-\cap U\subset \PH {\G}$ and
therefore $\PH{\g}^-\subset W$. However, every irreducible component of $W$
with non-empty boundary is contained in $\PH{\g}$ (cf. Proposition 2.3). 
\qed

\vfill\eject

\centerline{\bf 13. The Projective Spectrum.}\medskip

In this section we introduce a projective analogue of Gelfand's 
representation theorem for Banach algebras.  The relation of our 
construction to Gelfand's loosely  mirrors  the relation
of Grothendieck's  $\text{Proj}(R_*)$ of a graded ring $R_*$
to the spectrum $\text{Spec}(R)$ of an ordinary commutative ring $R$.

\Def {13.1} By a {\bf Banach graded algebra} we mean a graded normed
algebra 
$$
A_*\ =\ \bigoplus_{d\geq 0}A_d
$$
which is a direct sum of Banach spaces. Thus the norm on $A$ is a direct
sum $\|\bullet\|\  =\  \|\bullet\|_0  + \|\bullet\|_1  + \|\bullet\|_2 
+ \dots $ where $(A_d,\|\cdot\|_d)$ is complete, and 
$
\|a\cdot a'\|_{d+d'}\ \leq     \|a\|_d\ \|a'\|_{d'}
$
for all $a\in A_d$ and $a'\in A_{d'}$, or equivalently,
$$
\|a\cdot b\| \ \leq     \|a\|\ \|b\|
\qquad\text{for all }  a,b\in A_*.  
$$
A (degree-preserving) homomorphism of Banach graded algebras $\Psi:A_*\to
B_*$ is {\bf continuous} if there exists a constant $C>0$ such that
$||\Psi(a)|| \leq C^d ||a||$ for all $a\in A_d$ and all $d\geq0$.

\Ex{13.2}  Let $A_*=\bbc[t]$ be the algebra of polynomials in one variable
with 
$$
\|p(t)\|\ =\ \|\sum_{k=0}^n a_k t^k\| \ = \ \sum_{k=0}^n | a_k |.
$$
Note that the algebra automorphism $\bbc[t]$ determined by
$t\mapsto ct$, $c\neq 0$, is continuous with a continuous inverse.
This is the homogeneous coordinate ring of a projective
point.

\Ex{13.3} Consider a compact subset ${\X}\subset\bbp^n$ and let 
$A_d({\X}) = H^0(\bbp^n, \co{}(d))\bigl|_{\X}$ be the restriction of
holomorphic sections of $\co{\bbp^n}(d)$ to ${\X}$. Multiplication in
$A_*({\X})$  is induced by the   tensor product
$\co{}(d)\otimes\co{}(d')\to \co{}(d+d')$ and the  norm on $A_d$ is given 
by $$ \|\s\|_d \ =\  \sup_{\X}\|\s\|. $$

\redefine\X{X}

\Ex{13.4}  Let $\la\to X$ be a complex hermitian line bundle over a
locally compact topological space $X$  and let $A_d(X,\la) = \Gamma(X,
\la^d)$ denote the space of continuous sections of $\la^d$ with the
sup-norm. The Banach graded algebra $A_*(X, \la)$ will sometimes be called
the {\sl homogeneous coordinate ring of the polarized topological space }
$(X,\la)$.

\Def{13.5}  For a Banach graded algebra $A_*$ we denote  by 
$$H \equiv \Hom(A_*,\bbc[t])$$
the set of all continuous  degree-preserving  graded algebra homomorphisms 
$$
m: A_*  \arr\ \bbc[t].
$$
By definition of continuity, for each such $m$ there is a constant $C>0$
such that 
$$
|m(a)| \ \leq C^d \|a\|_d\qquad\qquad\text{for all } a\in A_d 
\text{ and all } d\geq0.
$$
We then set  $$H^\times = H-\{\underline0\}$$ where $\underline0$ denotes
the augmentation homomorphism $\underline0(a)=a_0$. (If we write $m\in H$
as $m=(m_0,m_1,m_2,...)$, then $\underline 0=(1,0,0,...)$.)

\Def{13.6} The  {\bf projective spectrum} of the graded algebra $A_*$ is
the quotient 
$$
\Proj(A_*)\ \equiv\ H^\times/\bbc^\times
$$ 
under the $\bbc^\times$-action on $H$ defined by 
$$
\phi_s(\{m_d\})\ =\ \{s^d m_d\}.
$$\medskip 
Given $m\in H^\times$ define
$$
\tbar m\tbar\ \equiv\ \inf\{ C \ :\ |m_d(a)| \leq C^d \|a\|_d \ \ \text{for
all } a\in A_d  \text{ and all } d>0\}.
$$
and set
$$
S(H)\ =\ \left\{m\in H^\times : \tbar m\tbar =1\right\}.
$$
The $\bbc^\times$-action
on $H^\times$ restricts to an $S^1$-action on $S(H)$. 
We introduce a topology on $\Proj(A_*)$ as follows. Embed
$$
S(H)\ \subset\ \prod_{d>0}\prod_{a\in A_d} D_a\ =\ D
\qquad\text{by } \ m\mapsto\{m_d(a)\}_{d,a}
\tag13.1$$
where $D_a=\{z\in\bbc:|z|\leq\|a\|\}$, and topologize $S(H)$ as a subspace
of $D$ with the product topology. The circle acts continuously on $S(H)$
by standard rotation in each factor, and    
$$
\Proj(A_*)=S(H)/S^1
$$
is given the quotient topology.  
Consider now the subset 
$$
B(H)\ =\ \left\{m\in H^\times : \tbar m\tbar \leq1\right\} \ \subset\ D
$$
embedded as in (13.1) above. 

\Prop{13.7} {\sl The set $B(H)\subset D$ is compact in the induced
topology. The quotient $\Proj(A_*)$ is compact if and only if $\underline0
\notin \overline{S(H)}$. 
}

\medskip

\Def{13.8} The algebra $A_*$ is called {\bf stable} if  $\underline0\notin
\overline{S(H)}$.  
\medskip

\noindent
{\bf Proof of Proposition 13.7.}
To see that $B(H)$ is closed in $D$ note that it is exactly the subset cut
out by the equations:    
$$\qquad\ \ \ \ 
z_{sa+ta'}=sz_a+tz_{a'} \qquad\text{for} \ \ a,a'\in A_d, d>0\ \ \text{and}
\ \ s,t\in \bbc
$$
$$\qquad\ \ 
z_a z_b \ =\ z_{ab}  \ \ \  \qquad\qquad\text{for} \ \ a\in A_d, \  b\in
A_{d'},\  d,d'>0.
$$
Evidently, $\Proj(A_*) = \pi(S(H)) = \pi(B(H)-\{\underline0\})$ where
$\pi:H^\times\to H^\times/\bbc^\times=\Proj(A_*)$ is the quotient map.
Hence, $\underline0\notin\overline{S(H)}$ implies that $\Proj(A_*) =
\pi(\overline{S(H)})$ is compact.  Conversely, if $\underline0
\in\overline{S(H)}$,
there is a net $Z_\a$ in $S(H)$ converging to $\underline0$.  If
$\Proj(A_*)$  were compact there would exist  a subnet $Z_\b$ with $\pi
Z_\b$ converging to some point $x\in \Proj(A_*)$. This however is
impossible, since the natural continuous map $D=\prod_a D_a\to\prod_a
[0,||a||]$, restricted to $S(H)$, descends to a continous map on
$\Proj(A_*)$. \qed

\medskip

The concept of stability is illuminated by considering the functions
$
||\bullet||_d : H^\times \ \arr\ \bbr^+
$
defined by
$$
||m||_d\ \equiv\ \sup\left\{\frac{|m_d(a)|}{||a||} : a\in A_d\right\}.
$$
These functions have the properties:
$$
||m||_{dd'} \ \geq\ ||m||_d ||m||_{d'}
$$
$$
|||m|||\ =\ \inf\{C: ||m||_d\leq C^d \text{ for all } d\}
$$
$$
|||m|||\ =\ 1\qquad\Leftrightarrow\qquad \sup_d ||m||_d = 1
$$
Finally note that 
$$
\underline0 \in \overline{S(H)} \ \ \Leftrightarrow\ \ \exists \text{ a
net }m^\a  \text{ in $S(H)$  s.t. } \lim_\a  m^\a_d(a)= 0 \text{ for
all } a\in A_d \text{ and  } d >0
$$
In particular,
$$
\exists \text{ a
net }m^\a  \text{ in $S(H)$  s.t. } \lim_\a || m^\a||_d= 0 \text{ for
all } d >0
 \ \ \Rightarrow\ \ \underline0 \in \overline{S(H)}.
$$

\vfill\eject

\def\cx{\Cal X}

\centerline{\bf 14. The Projective Gelfand Transform.}\medskip

Let $A_*$ be a Banach graded algebra and set $\cx = \Proj(A_*)$. Then
for each $d\geq 0$ there is  a   hermitian line bundle  
$$ 
\co{\cx}(d) \ \arr\ \cx
$$
associated to the principal $S^1$-bundle $S(H)\to\cx=S(H)/S^1$
by the character $t^d$ (considered
as a  homomorphism $S^1\to S^1=U(1)$).

Let $\ca(\cx,d)=\G(\cx,\co{\cx}(d))$  denote the space of continuous
sections of $\co{\cx}(d)$ equipped with the sup-norm topology.
Under tensor product, the direct sum
$$
\ca_*(\cx)\ \equiv\ \bigoplus_{d\geq0}\ca(\cx,d)
$$
becomes a Banach graded algebra.
 
Observe now that in terms of continuous functions on $S(H)$ we have
$$
\ca(\cx,d)\ =\ \{S:S(H)\to\bbc\ :\ S(\phi_t(m))=t^dS(m) \text {\ \
for all\ \ } t\in S^1\}.  
$$
Hence every element $a\in A_d$ gives rise to an element $\wh a\in
\ca(\cx,d)$ by setting
$$
\wh a(m)\ = \ m(a).
$$
This gives an embedding
$$
A_*\ \subset\ \ca_*(\cx)
\tag14.1$$
\Prop{14.1} {\sl  For all $a\in A_d$ one has 
$$
||\wh a||\ \leq \ ||a||.
$$
Thus the transformation (14.1) is a continuous injective
homomorphism of $A_*$ into the coordinate ring of the polarized
topological space $(\cx, \co{\cx}(1))$.
}
\pf Note that
$$
||\wh a||\ \equiv\ \sup_{[m]\in\cx} |\wh a(m)|
\ =\ \sup_{m\in S(H)} |\wh a(m)|
\ =\ \sup_{m\in S(H)} |m(a)|\ \leq\ \sup_{m\in S(H)} ||m||_d\cdot ||a||
\ \leq\ ||a|| 
$$
since $|||m|||=\sup_d ||m||_d=1$.\qed

\medskip\noindent
{\bf Question:} In the Gelfand case, one has
$||\wh a||=\lim_n||a^n||^{\frac1n}$. Is there an analogue here?

\vfill\eject

\centerline{\bf 15. Relation to the Projective Hull.}\medskip
Let ${\X}\subset \bbp^n$ be a compact subset and 
$A_*({\X})=\bigoplus_{d\geq0} H^0(\bbp^n, \co{\bbp^n}(d))\bigl|_{\X}$ the
algebra considered in 13.3. Set
$$
\cx\ \equiv\ \Proj(A_*({\X})).
$$
Note the natural embedding
$$
{\X}\ \hookrightarrow \ \cx
\tag15.1$$
which assignes to $x\in {\X}$ the   equivalence class of the multiplicative
functional $m_x:A_*({\X} )\to\bbc$ obtained by choosing an indentification
$\co x(1)\cong \bbc$ and setting $m_x(\cp)=\cp(x)$.

\Prop{15.1} {\sl The embedding (15.1) extends to a homeomorphism}
$$
\wh {\X}\ \cong\ \cx
$$
\pf Let $\pi:\bbc^{n+1}\to \bbp^n$ denote the projection and consider the
continuous mapping
$$
\pi^{-1}(\wh {\X})-\{0\}\ \arr\ H^\times
\tag15.2$$
given by $z\mapsto m_z$ where $m_z(p)=p(z)$ for homogeneous polynomials
$p$. This map is $\bbc^\times$-equivariant, i.e., $m_{tz} = t^dm_z$ on
$A_d$.  Note that 
$$
|||m_z|||\ =\ C([z])
$$
where $C=1/\rho$ is the best constant function (cf. Prop. 5.2). Therefore
the mapping (15.2) restricts to an $S^1$-equivariant mapping
$$
{\wh {\X}}_\rho\equiv \{z\in \pi^{-1}(\wh {\X} ): ||z||=\rho([z])\}\ \arr\
S(H). 
\tag15.3$$
which  induces a continuous mapping of the quotients
$$
\wh {\X}\ \arr\ \cx =  S(H)/S^1.
\tag15.4$$

A continuous inverse to this map is defined as follows.  For $m\in
H^\times$ consider the point 
$$
z\ =\ z_m\ = (m(Z_0),...,m(Z_n))
$$
where $Z_0,...,Z_n$ are the standard linear coordinates in $\bbc^{n+1}$.
Note that for any homogeneous polynomial $p(Z)\in \bbc[Z_0,...,Z_n]$ we
have
$$
m(p)\ =\ p(mZ_0,...mZ_n)\ =\ p(z)\ =\ m_z(p).
$$
Thus $m\mapsto z_m$ is a right inverse to (15.4). It is obviously also a
left inverse. \qed
\medskip

\Cor{15.2} {\sl The map (15.3) is an equivariant homeomorphism. In
particular, ${\X}$ is stable iff $A_*({\X})$ is stable.}\medskip

Note that in the   case considered here the Projective Gelfand
transformation simply maps the ``algebraic'' sections of $\co{}(d)\bigl|_{\X}$
(by extension) into the continuous sections of $\co{}(d)\bigl|_{\PH {\X}}$

\vfill\eject

\centerline{\bf 16. Finitely Generated Algebras.}\medskip
It is a classical fact that finitely generated Banach algebras 
correspond to polynomially convex subsets of $\bbc^n$.  We now show that 
analogously each finitely generated Banach graded algebra 
corresponds to a projectively convex subset  of $\bbp^n$.

\Prop{16.1}  {\sl
Let  $A_*$ be a Banach graded algebra generated by  elements
$a_0,...,a_N \in A_1$.  Then the algebra homomorphism
$$
\bbc[Z_0,...,Z_n]\ \arr\ A_*
\tag16.1$$
generated by  $Z_k\mapsto a_k$ is a continuous surjection which
induces a continuous injection 
$$
\Proj(A_*)\ \arr\ \bbp^n.
\tag16.2
$$
whose image is projectively convex (i.e., equal to its projective hull).}
\pf
By rescaling the generators (which induces a continuous isomorphism) we
may assume that $||a_k||=1$ for $k=0,...,n$. Observe now that
for any homogeneous polynomial $P=\sum c_\a Z^\a$,
$$
||P(a_0,...,a_n)||\ \leq\ 
\sum_\a |c_\a|\, ||a_0||^{\a_0}\dots||a_n||^{\a_n}
\ =\ \sum_{\a}|c_\a|\ \equiv ||P||_{\infty}
$$
and $||P||_{\infty}$ is equivalent to the standard  norm on
$\bbc[Z_0,...,Z_n]$ by Lemma A.1. The existence of the map (16.2) follows
immediately. It is induced by   the $\bbc^\times$-equivariant map
$$
H^\times\ \arr\ \bbc^{n+1}-\{0\}
$$
sending $m \mapsto (m(a_1),...,m(a_n))$. Under this map any homogeneous
polynomial $P(Z_0,...,Z_n)$ pulls back to
$P(ma_0,...,ma_n)
=m\{P(a_0,...,a_n)\}=\{P(a_0,...,a_n)\}^{\widehat{\ }}(m)$.  This is the
image of $P(a_0,...,a_n)$ in the homogeneous coordinate ring of  
$\Proj(A_*)$.

Let $\cx\subset \bbp^n$ denote the image of $\Proj (A_8)$. To see that 
$\PH{\cx} = \cx$ choose $[z]\in\PH{\cx}$.  By definition there is a
constant $C=C(z)$ such that 
$$
|P(z)|\ \leq  \ C^d\sup_{\cx}|P|\ =\ C^d\sup_{m\in S(H)}|mP(a)|\ \leq\ 
C^d||P(a)||
$$
 for all $P\in A_d$ and all $d$.  Hence $[z]\in \cx$.
\qed

\Note{16.2} The homomorphism (16.1) is only injective when $\Proj
(A_*)$ is Zariski dense.  In general we get a factoring of (16.1):
$$
\bbc[Z_0,...,Z_n]\ @>{\psi}>>\ {\wt A}_*\ @>{\phi}>>\ A_*
$$  
where ${\wt A}_*$ is the quotient algebra with the quotient 
norm in each degree, and where $\psi$ is an algebra isomorphism which is
continuous (but does not have a continuous inverse).  This induces
continuous injections: 
$$
\Proj(A_*)\ \arr\ \Proj( {\wt A}_*)\ \arr\ \bbp^n
$$
where $\Proj( {\wt A}_*)$ is the Zariski hull of $\cx=\Proj(A_*)$, i.e.,
the smallest algebraic subvariety containing $\cx$.
\medskip

\vfill\eject

\Note{16.3}  One can define a {\sl boundary} for $\cx=\Proj A_*\subset
\bbp^n$ to be a subset $\cx_0\subseteq\cx$ with $\PH{\cx_0}=\cx$.
As opposed to the affine case, there may be no unique minimal boundary.
For example if $\cx$ is an algebraic subvariety, then any open
subset, in fact any subset of positive $\omega$-capacity (cf [GZ]), is a
boundary. In particular boundaries can easily be disjoint.

On the other hand, for many of the examples considered here there is a
unique minimal boundary.

Note that by Theorem 12.8, the set $\cx-\cx_0$ for any boundary
$\cx_0$, satisfies the Local Maximum Modulus Principle for regular
functions.

\vskip .3in

\redefine\X{K}
\define\PSHo{\PSH_{\omega}}

\centerline{\bf 17. Projective Hulls on Algebraic Manifolds}\medskip

The projective hull of a subset can be defined abstractly in 
any projective variety.  Let $X$ be a compact complex manifold
provided with a  hermitian line bundle $\la$. 

\Def{17.1} The {\bf $\la$-hull } of a compact subset $\X\subset X$ is the
set $\PH {\X}_{\la}$ of all points $x\in X$ for which there exists a
constant $C=C_x$ such that $$
\|\s(x)\|\ \leq\ C^d\sup_{\X}\|\s\|
\tag17.1$$
for all $\s\in H^0(X,\co {} (\la^d))$ and all $d > 0$.
This set is independent of the metric on $\la$.
\medskip

Let $C_{\X,\la}:X\to(0,\infty]$ be the {\bf best constant}, defined at 
$x$ to be the smallest $C$ for which (17.1) holds, and set 
$$
{\V}_{\X,\la} = \log C_{\X,\la}.
$$  This function was studied by 
Guedj and Zeriahi [GZ] who introduced the following. Let $\omega$ denote
the curvature (1,1)-form of the hermitian connection on $\la$.

\Def{17.2} An upper semi-continuous function $v:X\to [-\infty,\infty)$ in
$L^1(X)$ is called  {\bf quasi-plurisubharmonic}\ \ if
$$
dd^cv + \omega\ \geq\ 0
\tag17.2$$
The convex set of such functions will be denoted $\PSHo(X)$
\medskip

Note that the   smooth  functions  $v\in\PSHo(X)$ are those
with the property that the hermitian metric $e^v\|\cdot\|$ has non-negative
curvature on $X$.

Note also that the u.s.c. function $\varphi=\frac 1 d \log \|\sigma\|$ with
$\sigma \in H^0(X, \co{}(\lambda))$ is in $\PSH_{\omega}(X)$ with
$dd^c\varphi+\omega =\frac 1 d \text{Div}(\sigma)$

\Theorem {17.3. [GZ]} {\sl Let $X, \la$ be as above  with $\la$ positive.
Then
$$
{\V}_{\X,\la}(x)\ =\ \sup\{v(x): v\in\PSHo(X)\ \ \text{ and }\ \
v\bigl|_{\X}\leq0\}. 
\tag17.3$$
Furthermore, the statements of Theorem 4.3 hold  with $\bbp^n$ replaced by
$X$.
}
\medskip

The $\lambda$-hull has the following elementary property.

\Lemma{17.4}   {\sl Let $\la\to X$ be a hermitian line bundle on a
compact complex manifold. Then for any compact set $\X\subset X$ and any 
$p\geq 1$
$$
\PH {\X}_{\la^p}\ =\ \PH {\X}_{\la} 
\and 
C_{\X,\la^p}\ =\ C_{\X, \la}^p
$$

}
\pf
It follow directly from the definitions that $\PH {\X}_{\la}\subseteq
\PH {\X}_{\la^p}$ and 
$C_{\X, \la^p}\ \leq C_{\X, \la}^p$. On the other hand, suppose
$x\in \PH {\X}_{\la^p}$ and $\s\in H^0(X,\la^d)$. Then
$$
\|\s(x)\|^p\ =\ \|\s^p(x)\|\leq C_{\X, \la^p}^d\sup_{\X}\|\s^p\|
\ =\ C_{\X, \la^p}^d \left(\sup_{\X}\|\s\|\right)^p,
$$
i.e., $\|\s(x)\|\leq C^d\sup_{\X} \|\s\|$ where
$C=(C_{\X,\la^p})^{\frac1p}$.
Hence, $x\in \PH{\X}_{\la}$ and 
$C_{\X, \la}\leq (C_{\X,\la^p})^{\frac1p}$.\qed
\medskip

We now examine the relationship of $\PH {\X}_{\la}$ to the projective
hull of $\X$ under projective embeddings related to $\la$.  
 Suppose $X\subset \bbp^N$ is
embedded   by the full space of sections of $\la$. Then 
for any Borel set
$\X\subset X$, one has that 
$
\PH {\X}_{X, \la}\ =\ \PH {\X},
$
and if $\la$ is given the metric induced from this embedding,
then
$
{\V}_{\X,X}\ =\ {\V}_{\X,\bbp^N}\bigl|_X.
$
(Of course ${\V}_{\X}\equiv \infty$ on $\bbp^N-X$.) 
This follows from the fact that any section of $\la^d$ 
is the restriction
of a section of $\co{\bbp^N}(d)$.
We now show that the hull remains unchanged if one embeds $X$ into
projective space by any subspace of $H^0(X,\co{}(\lambda))$.

\Prop{17.5}  {\sl  Suppose $X\subset \bbp^N$ is an embedding given
by a subspace  of sections of $\la$ for some $d$. Then for any compact
set $\X\subset X$ 
$$
\PH {\X}_{ \la}\ =\ \PH {\X}.
$$
Furthermore, is $\la$ is given the  metric induced from
this embedding, then
$$
{C}_{\X,\la}\ =\ {C}_{\X}\bigl|_X
$$
where $C_{\X}$ is the best constant function on $\bbp^N$.
(Of course ${C}_{\X}\equiv \infty$ on $\bbp^N-X$.)}

\Remark{17.6}   This result says essentially that the projective hull and
the associated extremal function of a subset $K\subset \bbp^n$ are 
intrinsic to any compact submanifold $X$ containing $K$.  

\medskip\noindent
{\bf Proof of Proposition 17.5}. 
This is a consequence of the following lemma.
We recall that $\la$ is {\bf very ample} if the sections of $\la$ give a
projective embedding of $X$.

\Lemma{17.7}  {\sl Let $X,\la$ be as above and suppose $f:Y\to X$ is a
holomorphic map from a compact complex manifold $Y$. Let $\mu = f^*\la$
with the induced metric.  Then:

(i)\ \ \ \  $f(\PH {\X}_{\mu}) \ \subseteq\ \PH {f(\X)}_{\la}$ \and
$f^*C_{f(\X), \la}  \ \leq\ C_{\X, \mu}$.
\smallskip

(ii)\ \ \  If $\la$ is very ample, then
$\PH {f(\X)}_{\la}\subseteq f(Y)$.

\medskip

(iii)\ \   If $\la$ is very ample and $f^*:H^0(X,\la^d)\to
H^0(Y,\mu^d)$ is surjective for all $d$, then }
$$
f(\PH {\X}_{\mu}) \ =\ \PH {f(\X)}_{\la}
\and
f^*C_{f(\X), \la}  \ =\ C_{\X, \mu}
$$
\pf
Suppose $y\in \PH {\X}_{\mu}$ and $\s \in H^0(X, \la^d)$. Then 
$$
\|\s(f(y))\|\  =\  \|(f^*\s)(y)\|\ \leq\ C_{\X,\mu}(y)^d \sup_{\X}
\|f^*\s\| \ =\ C_{\X,\mu}(y)^d \sup_{f(\X)} \|\s\|.
$$
Therefore, $f(y)\in \PH {f(\X)}_{\la}$ and 
$C_{f(\X), \la}(f(y)) \leq  C_{\X, \mu}(y)$. This proves (i).

For (ii) we note that 
$f(\X) \subset f(Y) \subset X\subset \bbp^N$ where the last embedding
is given by the sections of $\la$. By (i) we have 
$\PH {f(\X)}_{\la}\subset \PH {f(\X)}_{\co{\bbp^N}(1)}$ = the hull of
$f(\X)$ in $\bbp^N$. However, by Proposition 3.1(iii), the projective
hull is contained in the Zariski hull, and so 
$\PH {f(\X)}_{\la}\subseteq f(Y)$ as claimed.

For (iii) we suppose $x\in \PH {f(\X)}_{\la}$, so that
$
\|\s(x)\|\leq C_{f(\X),\la}^d\sup_{f(\X)}\|\s\|
$
 for all $\s\in H^0(X,\la^d)$ and all $d$.
Now by (ii), $x=f(y)$ for some $y\in Y$. Hence,
$$
\|(f^*\s)(y)\|\ =\ \|\s(fy)\|   \ \leq\
C_{f(\X),\la}(f(y))^d\sup_{f(\X)}\|\s\| \ =\ 
C_{f(\X),\la}(f(y))^d\sup_{\X}\|f^*\s\|
$$
for all $\s\in H^0(X,\la^d)$ and all $d$.  Therefore,
$$
\|\tau(y)\|\ \leq\   C_{f(\X),\la}(f(y))^d \sup_K\|\tau\|
\tag17.4$$
for all $\tau\in H^0(Y, \mu^d)$ and all $d$ by surjectivity.
Hence, $y\in \PH {\X}_{\mu}$ and so $x=f(y)\in f(\PH {\X}_{\mu})$.
Furthermore, by (17.4) we have $C_{\X, \mu}(y) \leq C_{f(\X),\la}(f(y))$.
Together with part (i) this completes the proof.\qed

\vfill\eject

\redefine\la{\Lambda}
\def\Do{F}

\centerline{\bf 18. Results for General K\"ahler Manifolds} \medskip

In this section we derive basic results concerning hulls of sets  in a
general setting. 

Let $X$ be a  K\"ahler manifold with  K\"ahler form $\omega$,
and fix a compact subset $K\subset X$.  Suppose $K\subset \Do\subset X$ with
$\Do$ compact and define
$$
\cs\ \equiv\ \PSHo(\Do)\ \equiv\ \{\vf \in C^\infty(X) : d d^c\vf + \omega
\geq 0 \text{ \ on }\ \Do\},
$$
the set of smooth functions on $X$ which are
quasi-plurisubharmonic   on $\Do$.

\Def {18.1}  For each $\la\geq0$ let ${\PH K_F}(\la)$
denote the set of all  $x\in \Do$ such that:
$$
\vf(x) \ \leq\ \sup_K \vf +\la \qquad\text{ for all}\ \ \vf \in \cs.
$$
The set 
$$
{\PH K_F}  \ = \ \bigcup_{\la\geq 0} {\PH K_F}(\la)
$$
will be called the {\bf $\omega$-quasi-plurisubharmonic hull  of $K$ in
$\Do$}.    When $X$ is compact we  set  
$\PH K(\la)=\PH K_X(\la)$ and $\PH K = \PH K_X$.

\medskip
Let $\cp_{1,1}(X)$ denote the set of positive currents of
bidimension (1,1) with compact support  on $X$, and let $\cm_K$
denote the set of probability measures on $K$.
%
%
%

\Theorem {18.2}  {\sl 
The following are equivalent.
\medskip

\qquad (A)\ \ \ $x \ \in\ {\PH K_F}(\la)$\medskip

\qquad (B) \ \ There exist $T\in \cp_{1,1}(X)$  with $M(T) \leq \la$ and
supp$(T)\subseteq F$ and a probablitiy 

\qquad \ \ \ \ \ \ \ \ measure $\mu \in \cm_K$ such that
$$
dd^c T \ =\ \mu - \delta_x
$$
}

For $\vf \in C^\infty(X)$ let 
$L_{\vf}$ denote the corresponding linear functional on
$\ce'(X)$.

\Lemma{18.3}  {\sl The following are equivalent.  \medskip

\qquad (i)\ \ \ $x \notin {\PH K_F}(\la)$

\medskip
\qquad (ii)\ \ \ There exists $\vf \in \cs$ with \ \  $\sup_K\vf +\la\ <\
\vf(x)$

\medskip
\qquad (iii)\ \ \ There exists $\vf \in \cs$ with \ \  $\int_K\vf\,d\mu +\la\
<\ \vf(x)$ for all $\mu\in\cm_K$.

\medskip
\qquad (iv)\ \ \ There exists $\vf \in \cs$ such that $\cm_K-\delta_x \
\subset \{L_\vf < - \la\}$
}

\pf
 We have (i) $\Leftrightarrow$ (ii) by definition. We have 
 (ii) $\Leftrightarrow$ (iii) because 
$$
\sup_K\vf = \sup_{\mu\in\cm_K}\int_K \vf\,d\mu.
$$
Note that $\cm_K$ is compact, so the strict inequality in (iii) implies
the strict inequality in (ii). Condition (iv) is just a restatement of
(iii).\qed\medskip

Consider the following subset of the compactly supported 0-dimensional
currents on $X$:
$$
C\ \equiv\ \{dd^cT\ :\ T \in \cp_{1,1}(X),\ M(T)\leq 1\text{ and } \supp
(T)\subseteq F \}.  $$
Obviously $C$ is a convex set containing the origin. It is easy to
see that $C$ is compact.

Recall that for a compact convex subset $\ck$
containing the origin   in a topological vector space $V$, the {\bf polar}
of $\ck$ is the  set $\ck^0\equiv \{L \in V^*: L\geq -1 \text{ on } \ck\}$.

\Prop{18.4}
$$
\cs\ =\ C^0
$$
\pf For $u\in C$ and $\vf \in \cs$ we have
$$
u(\vf) = (dd^cT)(\vf) = T(dd^c\vf+\omega) - T(\omega)\geq -T(\omega)\geq
-1. $$
Hence $\cs\subseteq C^0$.
Conversely suppose $\vf\in C^0$.  Then $-1\leq
u(\vf)=T(dd^c\vf+\omega)-T(\omega)$ or 
$$
0\leq T(dd^c\vf+\omega)+1-T(\omega)
$$
for all $T\in C$.  Taking $T=\delta_y \xi$ for $y\in F$ and $\xi$ a
positive simple unit (1,1)-vector at $y$,  we have $T(\omega)=M(T)=1$ and so
$(dd^c\vf+\omega)(\xi)\geq 0$. This proves that 
$(dd^c\vf+\omega)_y\geq 0$ for all $y\in F$. \qed
\medskip

Proposition 18.4 is equivalent to:
\Prop{18.4$'$}
$$
\vf\in \cs\ \ \Leftrightarrow\ \ \la C \subseteq \{ L_{\vf}\geq -\la\} 
$$
\pf 
$
\vf \in \cs  \ \Leftrightarrow\ \vf \in C^0 \ \Leftrightarrow\ 
C \subseteq \{u: L_{\vf}(u)\geq -1\} 
 \ \Leftrightarrow\ 
\la C \subseteq \{u: L_{\vf}(u)\geq -\la\}. 
$
\qed
\medskip

Note that
$$
\la C\ \equiv\ \{dd^cT \ :\ T \in \cp_{1,1}(X),\ M(T)\leq \la\text{ and }
\supp (T)\subseteq F \} 
$$
Combining Lemma 18.3 and the Proposition 18.4' yields:
\Prop{18.5}  {\sl The following are equivalent.  
\medskip

\qquad (i)\ \ \ $x \notin {\PH K_F}(\la)$

\medskip
\qquad (v)\ \ \   $\exists \vf \in C^\infty(X)$ with
$
\cm_K-\delta_x \ \subset \{L_\vf < - \la\}$ and $\la C \subset \{
L_{\vf}\geq -\la\}
$
}
\medskip

\noindent
{\bf Proof of Theorem 18.2.} The theorem can be restated as the
equivalence of: \medskip

\qquad (i)\ \ \ $x \notin {\PH K_F}(\la)$

\medskip
\qquad (vi)\ \ \  $\cm_K-\delta_x$ and $\la C$ are disjoint.
\medskip

\noindent
Obviously (v) $\Rightarrow$ (vi).  The Hahn-Banach Theorem states that 
(vi) $\Rightarrow$ (v). \qed

Suppose now that $X$ is compact and $F=X$, so that 
$\cs$ is the set of all quasi-plurisubharmonic functions on
$X$. In this case Theorem 18.2 can be strengthened so that 
$\supp(T)\subseteq {\PH K}^-$. This is the first main result of this
section.

\medskip

\Theorem {18.6}  {\sl Let $X$ be a compact K\"ahler manifold.  For any
compact subset $K\subset X$ the following are equivalent. \medskip

\qquad (A)\ \ \ $x \ \in\ \PH K(\la)$\medskip

\qquad (B) \ \ There exist $T\in \cp_{1,1}$ with $M(T)\leq \la$ and
$\supp(T)\subset {\PH K}^-$ such that
$$
dd^c T \ =\ \mu - \delta_x
$$
\qquad\qquad \ \ \ \ \ where  $\mu \in \cm_K$.
}

\pf 
That (B) implies (A) is already established in Theorem 18.2.

For the converse assume  $x \ \in\ \PH K(\la)$ but the equation in (B)  has
no solution.  Then by compactness there must exist  a compact subdomain 
$\Do$ with ${\PH K}^-\subset \Do^0$ such that there is no solution  $T\in
\cp_{1,1}(X)$ with $M(T)\leq \la$ and $\supp(T)\subseteq F$. Apply Theorem
18.2 to conclude that $x\notin {\PH K_F}(\la)$, that is,
there exists $\vf\in C^\infty(X)$ which is quasi-plurisubharmonic  on
$\Do$ with $\vf\leq 0$ on $K$ and $\vf(x)>\la$. 
It remains to find a function  $\wt\vf$ which is
quasi-plurisubharmonic on all of $X$ and
agrees with  $\vf$ on $\PH K(\la)$. Then $\wt\vf\leq0$
on $K$, and if  $x\in \PH K(\la)$, then $\wt\vf(x)=\vf(x)>\la$, which is a
contradiction. 

\Prop {18.7}  {\sl
Assume $X$ is a compact K\"ahler manifold and $\la>0$.  Suppose $\vf\in
C^\infty(X)$ is  quasi-plurisubharmonic on a
neighborhood   of ${\PH K}^-$.  Then there exists a $C^\infty$
quasi-plurisubharmonic function $\wt \vf$ on $X$ which agrees with $\vf$
in a neighborhood of $\PH K(\la)$.
}

\Lemma {18.8} {\sl
Assume $X$ is a compact K\"ahler manifold and $\la>0$.
For each open neighborhood $U$ of ${\PH K}^-$ and each $N$ large, there
exists a $C^\infty$ quasi-plurisubharmonic function $\psi$ on $X$ with
$\psi>N$ on $X-U$ and $\psi < -N$ on some neighborhood of
$\PH K(\la)$.
}

\medskip

\noindent
{\bf Proof.}  Note that if $\vf\in \PSHo(X)$  and $\vf\leq0$
on $K$, then $\vf \leq \la$ on  $\PH K (\la)$.  For each $y\in X-U$, since
$y\notin {\PH K}^-$, there exists 
$\psi\in \PSHo(X)$ with $\psi\leq0$ on $K$ and $\psi(y) > 2N+\la$.  Set
$V_y=\{x\in X : \psi(x)>2N+\la\}$. Extract a finite subcover $V_1,...,V_r$
of $X-U$ with associated functions  $\psi_1,...,\psi_r$.  Let $\psi =
\max\{\psi_1,...,\psi_r\} \in\PSHo(X)$ (see [GZ, Prop. 1.3]).    Then
$\psi > 2N + \la$ on a neighborhood of $X-U$ and $\psi\leq 0$ on $K$. 
Therefore $\wt{\psi} = \psi-N-\la$ satisfies $\wt{\psi}\leq -N $ on $\PH
K(\la)$ and $\wt{\psi}>N$ on a neighborhood of $X-U$. Finally replace
$\wt{\psi}$ by $\psi= \wt{\psi}-\delta$ with $\delta>0$ sufficiently 
small that we still have   ${\psi}>N$ on a neighborhood of $X-U$.
Then $\psi<-N$ on some neighbohood of $\PH K(\la)$.\qed

\medskip
\noindent
{\bf Proof of Proposition 18.7.} Suppose $\phi$ is
quasi-plurisubharmonic on  $U \supset {\PH K}^-$.
Now pick $N$ so that $|\vf| < N$
on $\overline U$.  Then $\wt {\vf} \equiv \max\{\vf , \psi\}$ satisfies:
\smallskip

\qquad 1) \ \ $\wt{\vf} = \vf$ \ \ in a neighborhood of $\PH K(\la)$,
\smallskip

\qquad 2) \ \ $\wt{\vf} = \psi$ \ \ in a neighborhood of $X-U$. 
\hfill\blbx
\!\!\!\!\!\!\!\!\!\!\!\!\!\!\!\!
\!\!\!\!\!\!\!\!\!\!\!\!\!\!\!\!\!\!\!\!\!\!\!\!\!\!\!\!\!     
\blbx

\Remark{} The proofs of Proposition 18.7 and Lemma 18.8 only produced a
continuous function since in general $\max\{\vf,\psi\}$ is only
continuous. However, 
$$
\max\{\vf,\psi\}\ =\ \lim_{n\to\infty}\frac
1n\log \left( e^{n\phi}+e^{n\phi} \right) 
$$
can be approximated by smooth quasi-plurisubharmonic functions
(See [GZ], [D$_2$]).

\medskip

Theorem 18.6 can be extended to the non-compact case.
On any $X$ we continue to define  $\PH K$ and $\PH K(\la)$   as in 18.1
with $\Do=X$.

\eject

\Theorem{18.9} {\sl   Let $X$ be a non-compact K\"ahler manifold.      
Then for any compact subset $K\subset X$ with  $\PH K\subset\subset X$ the
following are equivalent.  \medskip

\qquad (A)\ \ \ $x \ \in\ \PH K(\la)$\medskip

\qquad (B) \ \ There exist $T\in \cp_{1,1}(X)$ with $M(T)\leq \la$ and
$\supp(T)\subset {\PH K}^-$ such that
$$
dd^c T \ =\ \mu - \delta_x
$$
\qquad\qquad \ \ \ \ \ where $\mu \in \cm_K$.
}

\pf
Suppose $T$ is the current asserted in (B) and choose $\vf\in\PSHo (X)$.
Then $\int \vf\, d\mu-\vf(x)= d d^cT(\vf)=T(d d^c \vf)=T(d d^c \vf+\omega)
-T(\omega)\geq -\la$, and so $x\in \PH K(\la)$.

Suppose now that  $x\in \PH K(\la)$ and (B) does not hold. Then there
must  exist  a compact subdomain $\Do$ with ${\PH K}^-\subset \Do^0$ such
that there exists no solution $T\in \cp_{1,1}(X)$ with $M(T)\leq \la$ and
$\supp(T)\subseteq F$.
 Hence, by Theorem
18.2 there exists $\vf\in\PSHo(\Omega)$ with $\vf\leq0$ on $K$ and
$\vf(x)>\la$. Choose a larger compact subdomain $D$ with $\Do\subset 
\subset D^0$.  Fix $N>\sup_{\Do}|\vf|$. The argument given for Lemma 18.8
shows that there exists $\psi\in \PSHo(X)$ with $\psi<-N$ on a
neighborhood of $\PH K(\la)$ and $\psi>N$ on a neighborhood of $D-\Do^0$.
Define $\widetilde{\vf}\in \PSHo(X)$ by
$$
\widetilde{\vf}=\cases
\max\{\vf,\psi\} \qquad\text{ on }D \\
\ \ \ \ \ \ \psi\  \qquad \qquad\text{ on } X-D
\endcases
$$
and note that $\widetilde{\vf}=\vf$ in a neighborhood of $\PH K(\la)$.
However, $\widetilde{\vf}\leq0$ on $K$ and $\widetilde{\vf}(x)>\la$, so
$x\notin \PH K(\la)$, a contradiction.
\qed

\Remark{18.10} {\sl The analogues of Corollaries 11.3 and 11.4, and
Theorems 12.1 and 12.3 hold in this context. Moreover, the following
analogue of Theorem 12.8 holds.\medskip

Fix ${\X}^{\text{compact}}\subset X$ and $U^{\text{open}}\subset\subset
\Omega^{\text{open}}\subset X$ where $\Omega$ is analytically equivalent
to a Runge domain in $\bbc^n$. Then}
$$\boxed{
{\PH {\X}}^- \cap U \ \subseteq \text{$\Omega$-HolomorphicHull}
\left\{({\PH {\X}}^-\cap\partial
U)\cup ({\X}\cap U)\right\}
}
 $$
{\sl with equality if $C_{\X}$ is bounded on $\PH {\X} \cap \overline U$.}

\Remark{18.11}  Much of the discussion of section 4 holds in this
general context. The capacity of Dinh-Sibony
[DiS] was introduced for any K\"ahler manifold $X$ and  the Theorem 4.3 of
Guedj-Zeriahi holds there. Furthermore, Dinh-Sibony [DiS] proved that 
any analytic subvariety $Z\subset X$ is always globally
$\omega$-pluripolar. Hence, if $K\subset Z$, then $\PH K \subset Z$, and
so  $\PH K$ is contained in the ``analytic hull'' of $K$, that is,  the
intersection of all subvarieties of $X$ which contain $K$.

\redefine\la{\lambda}

\vfill\eject

\centerline{\bf Appendix A.\ \  Norms on $A_*(\bbp^n)$.}\medskip
From one point of view this paper is simply concerned with the study 
of  equivalence classes of  norms on the graded algebra
$\bbc[Z_0,...,Z_n]$. Two norms $||\bullet||$ and $||\bullet||' $ are
equivalent if there exists a constant $C>0$ such that 
$
\frac 1 {C^d} ||a|| \leq ||a||' \leq C^d ||a|| 
$
for all $a \in A_d$ or equivalently, if the identity map
$(A_*,||\bullet||)\to (A_*,||\bullet||')$ is continuous in both directions.
 There are many norms
equivalent to the standard one given by (2.7).  We examine some of them
here.

Let $\Omega\subset \bbc^n$ be a closed bounded convex set and define
$$
||P||_{\Omega}\ \equiv\ \sup_{\Omega}|P|
$$
for $P\in\bbc[Z_0,...,Z_n]_d$. Then obviously 
$$\aligned
(i)\ \ \ &\Omega_1\subset \Omega_2 \ \ \Rightarrow 
||P||_{\Omega_1} \ \leq\  ||P||_{\Omega_2}  \\
(ii)\ \ \ &||P||_{t\Omega} \ =\ |t|^d||P||_{\Omega} 
\endaligned
$$
for all $P\in\bbc[Z_0,...,Z_n]_d$. It follows easily that {\sl all these
norms are equivalent}.

Another interesting norm on $\bbc[Z_0,...,Z_n]$ is defined 
on $P(Z)=\sum_{|\a|=d} c_{\a}Z^{\a}$ by
$$
||P||_{\infty}\ \equiv\ \sum_{|\a|=d} |c_{\a}|
$$
\Lemma{A.1} {\sl The norms $||\bullet||_{\infty}$ and
$||\bullet||_{\Omega}$ are equivalent.}
\pf We shall work with the polydisk $\Omega = \{Z\in\bbc^{n+1}:|Z_k|\leq1
\text{ for all } k\}$.  Note that 
$$\aligned
||P||_{\Omega} \ &= \ \sup_{|Z_0|=\dots=|Z_n|=1}|P(Z)|
\ =\ \sup_{|Z_0|=\dots=|Z_n|=1}\bigl|\sum_{|\a|=d} c_{\a}Z^{\a}\bigr|
  \\
 & \ \leq \sup_{|Z_0|=\dots=|Z_n|=1}\sum_{|\a|=d}|c_{\a}|\,|Z^{\a}|
\ =\ \sum_{|\a|=d}|c_{\a}|\ =\ ||P||_{\infty}.
\endaligned
$$
For the converse assume inductively that
$$
||P||_{\infty}\ \leq\ C^d||P||_{\Omega} 
$$
for all $P\in\bbc[Z]_d$, and all $ d\leq N-1$ where $C=(n+1)4^{n+1}$.
Fix $P=\sum_{\a} c_{\a}Z^{\a}\in\bbc[Z]_N$ and note that
$$
\frac{\partial P}{\partial Z_k}\ =\ \sum_{|\a|=N} \a_kc_{\a}Z^{\a-\e_k}\
\in\ \bbc[Z]_{N-1}. $$
Hence by induction
$$
\sum_{\a} \a_k|c_{\a}|\ \leq \ C^{N-1}\left|\left|\frac{\partial P}{\partial
Z_k}\right|\right|_{\Omega}. $$
Now for $Z\in \Omega$,
$$
\frac{\partial P}{\partial Z_k}(Z)\ =\ \left(\frac1{2\pi
i}\right)^{n+1}\underset{|\zeta_0|=2}\to {\int}\dots
\underset{|\zeta_n|=2}\to {\int} \frac{P(\zeta)\,d\zeta_0\dots d\zeta_n}
{(\zeta_k-z_k)\prod_{j=0}^n(\zeta_j-z_j)}
$$
from which it follows that 
$$\aligned
\left|\left|\frac{\partial P}{\partial Z_k}\right|\right|_{\Omega}
\ &\leq\ \left(\frac1{2\pi}\right)^{n+1}\int_0^{2\pi}\dots\int_0^{2\pi}
|P(2e^{i\theta_0},...,2e^{i\theta_n})|\,2^{n+1}d\theta_0\dots d\theta_n\\
&\ \\
&\leq 2^{n+1}\sup_{|\zeta_0|=\dots=|\zeta_n|=2} |P(\zeta)|
\leq 4^{n+1}\sup_{|\zeta_0|=\dots=|\zeta_n|=1} |P(\zeta)|
= 4^{n+1}||P||_{\Omega}
\endaligned
$$
Therefore we have
$$\aligned
||P||_{\infty}\ &=\ \sum_\a |c_\a|\ \leq \sum_{k=0}^n\a_k|c_\a|
\ =\ \sum_{k=0}^n
\left|\left|\frac{\partial P}{\partial Z_k}\right|\right|_{\infty}
\\
&\leq\ C^{N-1}\sum_{k=0}^n
\left|\left|\frac{\partial P}{\partial Z_k}\right|\right|_{\Omega}
\ \leq \ C^{N-1} (n+1)4^{n+1}||P||_{\Omega}\ =\ C^N||P||_{\Omega}.
\endaligned
$$
as desired. \qed

\centerline{\bf References}

\vskip .2in

\ref \key A$_1$ 
\by  \ \ H. Alexander \paper Polynomial approximation and
hulls in sets of finite linear measure in $\bbc^n$
  \jour Amer. J. Math.
 \vol 93    \yr 1971    \pages  65-74 \endref

\smallskip

\ref \key A$_2$  \by \ \ H. Alexander \paper Projective capacity.
\inbook Recent developments in several complex variables (Proc. Conf.,
Princeton Univ., 1979)\pages 3-27 
\publ  Ann. of Math. Studies No. 100,  Princeton Univ. Press, Princeton,
N.J. \yr 1981  \endref

\smallskip

\ref \key  AW$_1$  \by \ \ \ \ \  H. Alexander and J. Wermer \book Several Complex
Variables and Banach Algebras
 \publ Springer-Verlag \publaddr New York   \yr  1998 \endref

 \smallskip

\ref \key  AW$_2$  \by \ \ \  \  H. Alexander and J. Wermer \paper Linking numbers
and boundaries of varieties
  \jour Ann. of Math.
 \vol 151  \yr  2000  \pages 125-150\endref

 \smallskip

\ref \key  BT  \by \ \ E. Bedford and B. A. Taylor \paper A new capacity for
plurisubharmonic functions
 \jour Acta Math.
 \vol 149   \yr  1982  \pages 1-40 \endref

 \smallskip



\ref \key  D$_1$  \by \ \ J.-P. Demailly \paper Estimations $L^2$ pour
l'op\'erateur $\overline{\partial}$ d'un fibr\'e vectoriel holomorphe
semi-positif au-dessus d'une vari\'et\'e k\"ahl\'erienne compl\`ete
 \jour Ann. Sci. E. N. S. (4)
 \vol 15   \yr  1982  \pages 457-511 \endref

 \smallskip

\ref \key  D$_2$  \by \ \   J.-P. Demailly \paper Regularization of closed
positive currents and intersection theory
 \jour  J. Algebraic Geom. 
 \vol 1   \yr  1992  \pages 361-409 \endref

 \smallskip

\ref \key  DPS  \by \ \ \ \ J.-P. Demailly, T. Peternell and M Schneider \paper
Pseudo-effective line bundles on compact K\'ahler manifolds
 \jour Internat. J. Math.
 \vol 12   \yr  2001 \pages 689-741 \endref

 \smallskip

\ref \key  DF  \by \ \ \ K. Diederich and J. E. Forna\!ess \paper
A smooth curve in $\bbc^2$ which is not a pluripolar set
 \jour Duke Math. J.  \vol  49 \yr 1982  \pages 931-936
\endref

 \smallskip

\ref \key  DL  \by \ \ T.-C. Dinh and M. Lawrence \paper
 Polynomial hulls and positive currents
 \jour Ann. Fac. Sci de Toulouse \vol  12 \yr 2003   \pages 317-334
\endref

 \smallskip

\ref \key  DiS  \by \ \ \ T.-C. Dinh and N. Sibony \paper
 Distribution des valeurs des transformations meromorphes et applications
 \jour  Preprint ArXiv:math. DS/0306095 \yr 2003    
\endref

 \smallskip

\ref \key  Do  \by \ \ P. Dolbeault \paper On holomorphic chains with given 
boundary in $\bbc\bbp^n$
 \jour Springer  Lecture Notes, no. 1089,
   \yr 1983   \pages 1135-1140\endref

 \smallskip

\ref \key  DH$_1$  \by  \ \ \ \ P. Dolbeault and G. Henkin \paper
Surfaces de Riemann de bord donn\'e dans $\bbc\bbp^n$,
\inbook Contributions to complex analysis and analytic geometry
 \jour Aspects of Math.
\publ  Vieweg \vol 26 \yr 1994   \pages 163-187 \endref

\smallskip

\ref \key  DH$_2$  \by \ \ \  \ P. Dolbeault and G. Henkin \paper
 Cha\^ines holomorphes de bord donn\'e dans $\bbc\bbp^n$
 \jour Bull.  Soc. Math. de France \vol  125 \yr 1997   \pages 383-445
\endref

 \smallskip

\ref \key  DS  \by \ \ \ J. Duval and N. Sibony \paper
Polynomial convexity, rational convexity and currents
 \jour Duke Math. J. \vol  79 \yr 1995   \pages 487-513
\endref

 \smallskip

\ref \key  Fa$_1$  \by \ \ \ \ B. Fabr\'e \paper
Sur l'intersection d'une surface de Riemann avec des hypersurfaces  
alg\'ebriques
  \jour C. R. Acad. Sci. Paris
 \vol  322 S\'erie I   \yr 1996    \pages 371-376   \endref

\smallskip

\ref \key  Fa$_2$  \by \ \ \  B. Fabr\'e \paper
On the support on complete intersection 0-cycles
  \jour The Journal of Geometric Analysis
 \vol  12    \yr 2002    \pages 601-614   \endref

\smallskip

\ref\key  F \by  \   H. Federer\book Geometric Measure 
Theory\publ 
 Springer--Verlag\publaddr New York \yr 1969\endref

 \smallskip

\ref \key  G  \by \   I. M. Gelfand \paper
Normierte Ringe
 \jour Math. Sb.  (N.S.) \vol 9 (51)  \yr 1941 \pages   3-24  
\endref

 \smallskip

\ref \key  GZ  \by \ \ \ V. Guedj and A. Zeriahi \paper
Intrinsic capacities on compact K\"ahler manifolds
 \jour Preprint Univ. de Toulouse   \yr 2003  
\endref



%

\smallskip

\ref \key  HL$_1$  \by \ \ \ \  F.R. Harvey and H.B. Lawson, Jr. \paper
On boundaries of complex analytic varieties, I  \jour Ann. of Math.
 \vol 102 \yr 1975  \pages 223-290\endref

\smallskip

\ref \key  HL$_2$  \by \ \ \ \ F.R. Harvey and H.B. Lawson, Jr. \paper
On boundaries of complex analytic varieties, II  \jour Ann. of Math.
 \vol 106 \yr 1977  \pages 213-238\endref

\smallskip

\ref \key  HL$_3$  \by \ \ \  F.R. Harvey and H.B. Lawson, Jr. \paper
Boundaries of varieties in projective manifolds  \jour Jour. of Geom.
Analysis \vol 14 no. 4 \yr 2005\endref

\smallskip

\ref \key   HL$_4$  \by \ \ \  F.R. Harvey and H.B. Lawson, Jr. \paper
Projective hulls, projective linking and boundaries of
positive holomorphic chains in projective manifolds  \jour Stony
Brook preprint \yr 2005\endref

\smallskip

\ref\key  Hi  \by  \ \  E. Hille \book Analytic Function Theory, Vol II
\publ Ginn and Co.
 \publaddr Boston \yr 1962 \endref

\smallskip

\ref\key  Ho \by  \  L. H\"ormander \book An Introduction to Complex
Analysis in Several  Variables
\publ Van Nostrand
 \publaddr Princeton, N. J.  \yr 1966 \endref

\ref \key  LMP  \by \ \ \ \  N. Levenberg, G. Martin and E. A. Poletsky \paper
Analytic disks and pluripolar sets 
 \jour  Indiana Univ. Math. J. \vol  41 \yr 1992   \pages 515-532
\endref

 \smallskip

\ref \key  R  \by \  M. Rosenlicht \paper
Generalized Jacobian Varieties
 \jour Ann. of Math. \vol  59  \yr 1954   \pages 505-530
\endref
\smallskip

\ref \key  S  \by \ A. Sadullaev \paper
 An estimate for polynomials on analytic sets
 \jour Math. USSR Izvestia \vol  20 No. 3 \yr 1983   \pages 493-502
\endref

 \smallskip

\ref \key  Sib  \by \ \ N. Sibony  \paper multi-dimensional analytic
structure in the spectrum of a uniform algebra \inbook 
Spaces of analytic functions. (Seminar on Functional Analysis and Function
Theory, Kristiansand, 1975) \publ Springer Lect. Notes in Math. No. 512
   \yr 1976 \pages 139-165 \endref

 \smallskip
\ref \key  SW  \by \ \  N. Sibony and P.-M. Wong \paper Some results on
global analytic sets
 \jour Springer S\'eminaire Pierre Lelong-Henri Skoda (Analyse). Ann\'ees
1978-79. pp. 221-237,  Lect. Notes in Math.
 \vol 822   \yr  1980    \endref

 \smallskip

\ref \key  Si  \by \   J. Siciak \paper On some extremal functions and
their applications in the theory of analytic functions in several
variables
 \jour Trans. Amer. Math. Soc.
 \vol 105   \yr  1962  \pages 322-357 \endref

 \smallskip

\ref \key  Wie  \by\ \ \  J. Wiegerinck \paper Pluripolar sets: hulls and
completeness
 \jour Actes des Rencontres d'Analyse Complexe (Poitiers-Futuroscope,
1999) \publ Atlantique, Poitiers,
      2002  \pages 209-219 \endref

 \smallskip

\ref \key  W$_1$  \by \ \ \  J. Wermer \paper The hull of a curve in $\bbc^n$
 \jour Ann. of Math.
 \vol 68   \yr  1958  \pages 550-561 \endref

 \smallskip

\ref \key  W$_2$  \by \ \ \,  J. Wermer \paper The argument principle and
boundaries of analytic varieties
 \jour Operator Theory: Advances and  Applications
   \vol 127   \yr  2001  \pages 639-659 \endref

 \smallskip

\ref \key  Z  \by \ \  A. Zeriahi \paper A criterion of algebraicity  for
Lelong classes and analytic sets
 \jour Acta Math.
 \vol 184   \yr  2000  \pages 113-143 \endref

 \smallskip

\end